\documentclass[11pt]{article}

\usepackage{a4}
\usepackage{geometry}
\usepackage{amsfonts}
\usepackage{graphicx}
\usepackage{epstopdf}
\usepackage{amsmath}
\usepackage{amsthm}
\usepackage{stmaryrd}
\usepackage{bm}
\usepackage{hyperref}
\usepackage{pgfplots}

\ifpdf
  \DeclareGraphicsExtensions{.eps,.pdf,.png,.jpg}
\else
  \DeclareGraphicsExtensions{.eps}
\fi

 \newtheorem{theorem}{Theorem}[section]

 \theoremstyle{definition}
 
 \theoremstyle{remark}
 \newtheorem{remark}[theorem]{Remark}
 \numberwithin{equation}{section}

\title{Upscaling of a Cahn--Hilliard Navier--Stokes Model with Precipitation and Dissolution in a Thin Strip\thanks{Acknowledgment:
Funded by the Deutsche Forschungsgemeinschaft (DFG, German
Research Foundation) -- Project Number 327154368 -- SFB 1313, the Research Foundation Flanders (FWO) -- Project G0G1316N, and the Hasselt University -- Project BOF19BL12.}}

\author{Lars von Wolff
\thanks{Institute of Applied Analysis and Numerical Simulation,
		University of Stuttgart,
		Pfaffenwaldring 57,
		70569 Stuttgart, 
		Germany
  (\mbox{lars.von-wolff@mathematik.uni-stuttgart.de}).
 }
\and Iuliu Sorin Pop
\thanks{Faculty of Sciences,
		Hasselt University,
		Agoralaan Gebouw D,
		3590 Diepenbeek,
		Belgium
  (\mbox{sorin.pop@uhasselt.be}).
 }
}

\usepackage{amsopn}

\usepackage{tikz}
\usepackage{cases}
\usetikzlibrary{arrows.meta}  
\usetikzlibrary{intersections}

\newcommand{\NN}{{\mathbb{N}}}         
       
\newcommand{\RR}{{\mathbb{R}}}

\newcommand{\vv}{{\mathbf{v}}}
\newcommand{\Phiv}{{\mathbf{\Phi}}}
\newcommand{\nv}{{\mathbf{n}}}
\newcommand{\Jv}{{\mathbf{J}}}
\newcommand{\Sv}{{\mathbf{S}}}

\newcommand{\ev}{{\mathbf{e}}}
\newcommand{\xv}{{\mathbf{x}}}

\newcommand{\tauv}{{\bm{\tau}}}

\newcommand{\diff}{\,d}											
\newcommand{\eps}{{\varepsilon}}								
\newcommand{\set}[1]{\left\{#1\right\}}            				

\newcommand{\jump}[1]{\llbracket #1 \rrbracket}

\renewcommand{\r}[1]{{#1}_{\text{ref}}}
\newcommand{\overbar}[1]{\mkern 1.5mu\overline{\mkern-1.5mu#1\mkern-1.5mu}\mkern 1.5mu}

\hypersetup{
  pdftitle={Upscaling of a Cahn--Hilliard Navier--Stokes Model with Precipitation and Dissolution in a Thin Strip},
  pdfauthor={L. von Wolff, I. S. Pop},
  pdfkeywords={Two-phase flow; Fluid flow with reactive transport; Precipitation/dissolution;  Phase-field models; Upscaled models; Asymptotic analysis}
  pdfsubject={35R35, 35Q35, 76D05, 35C20, 35B27, 76M50}
}

\begin{document}
\maketitle

\begin{abstract}
We consider a phase-field model for the incompressible flow of two immiscible fluids. This model extends widespread models for two fluid phases by including a third, solid phase, which can evolve due to e.g. precipitation and dissolution.

We consider a simple, two-dimensional geometry of a thin strip, which can still be seen as the representation of a single pore throat in a porous medium. Under moderate assumptions on the P\'eclet number and the capillary number, we investigate the limit case when the ratio between the width and the length of the strip is going to zero. In this way and employing transversal averaging,  we derive an upscaled model. The result is a multi-scale model consisting of the upscaled equations for the total flux and the ion transport, while the phase-field equation has to be solved in cell-problems at the pore scale to determine the position of interfaces. We also investigate the sharp-interface limit of the multi-scale model, in which the phase-field parameter approaches 0. The resulting sharp-interface model consists only of Darcy-scale equations, as the cell-problems can be solved explicitly. Notably we find asymptotic consistency, that is the upscaling process and the sharp-interface limit commute. We use numerical results to investigate the validity of the upscaling when discontinuities are formed in the upscaled model.
\end{abstract}
\providecommand{\keywords}[1]{{\textit{Key words:}} #1\\ \\}
\providecommand{\class}[1]{{\textit{AMS subject classifications:}} #1}
\class{35R35, 35Q35, 76D05, 35C20, 35B27, 76M50 }

\section{Introduction}

Multi-phase flow and reactive transport in porous media are encountered in many important fields, including geological $CO_2$ sequestration, geothermal energy, groundwater management, oil recovery and ion exchange in fuel cells. While the modelling of multi-phase flow is itself a challenging task, the examples given before have in common that the solid matrix of the porous medium can change in time due to processes like precipitation or dissolution, which, in turn influence the flow behaviour. 

Another common point of the processes mentioned before is that they are taking place in a porous medium. In this case, two different length scales are encountered. At the scale of pores, each phase (solid, or fluid) is identified clearly, occupying certain positions in well defined volumes. At the Darcy-scale, which is often the scale of main interest, averaged quantities are used to describe the behaviour of the system . 

In detail, we are interested here in the situation where two immiscible fluid phases are occupying the pore space of a porous medium. One fluid phase contains ions that can precipitate at the fluid-solid interfaces. This leads to the formation of a precipitate layer at the pore walls, which reduces the space available for the fluid. The reverse process, that is the dissolution of the mineral phase into the fluid phase, is also allowed. In this case, the volume of the precipitate is reduced, while the volume available for flow is increased, and more ions are dissolved in the fluid phase.

To model this process at the pore scale, one uses the conservation of mass, momentum and of the dissolved ions in each phase. Since the spaces occupied by each of the two fluids, and of the mineral as well can change over time, two free boundaries are encountered at the pore scale. These free boundaries are separating the different phases. 

Different approaches have been proposed for developing the corresponding mathematical models. For a simple geometry, which is basically a long, thin strip (in two spatial dimensions) or tube (in three dimensions) the free boundaries can be viewed as functions of one or two variables. In this sense we mention \cite{Tycho1} for a model describing precipitation and dissolution but for one fluid phase, which has been extended in \cite{Agosti, Carina15, KNP, KundanALE}, and  \cite{AndroTube, AndroStrip, Lunowa, Picchi, Sohely} for unsaturated single-phase flow or two-phase flow models. 

For more complex geometries, level sets can be employed to describe the evolution of the free boundaries. In this respect we refer to \cite{noorden2009b}, as well as to \cite{Carina16, Schulz1, Schulz2}, all considering models for precipitation and dissolution in a water-saturated porous medium. 

When applying any of both approaches mentioned before, one has to deal with (freely) moving interfaces. This makes not only the mathematical analysis, but also the development of efficient numerical scheme a challenging task. 
Alternatively, one can use phase-fields to approximate the interfaces between phases by diffuse transition zones with small positive width. The phase-fields are smooth approximations of the indicator function of each phase. The evolution of the phase fields is usually derived as the gradient-flow to a free energy and, in the limit case when passing the diffuse interface parameter, one should recover the original, free boundary model. 

Commonly used phase-field models are involving either the Allen--Cahn equation \cite{allencahn} or the Cahn--Hilliard equation \cite{cahnhilliard}. While the Allen-Cahn equation is of second order and ensures that the phase-field indicators remain essentially bonded by zero and one, it is not conservative. Therefore here we focus on the Cahn--Hilliard equation, which is of fourth order but conservative for the phase-field indicators. 

Models coupling the Cahn--Hilliard equations and the incompressible Navier--Stokes equations have been developed for two fluid phases \cite{agg}, three fluid phases \cite{Boyer06, Boyer10}, and more than three fluid phases \cite{Boyer14, stinner19}. For the description of fluid-solid interfaces, the Navier--Stokes equations can be solved in the fluid volume fraction and a velocity of zero is assigned to the solid phase \cite{Beckermann, Beckermann04}. Phase-field models are also used in \cite{Banas, Bunoiu, Daly, Metzger, Schmuck1, Schmuck2} as pore-scale models for two-phase flow in porous media, and further Darcy-scale models are derived. Kinetic reactions at phase boundaries have been introduced in \cite{noorden2011, magnus}. The pore-scale model in \cite{magnus} includes two immiscible fluid phases and a mineral one, but the fluid phases only move due to curvature effects. Also, the corresponding Darcy-scale model is derived by homogenization techniques. More recently, phase-field models that couple precipitation and dissolution with fluid flow have been developed in \cite{CarinaLarsSorin} (for one fluid phase, and for which the Darcy-scale model is derived), and \cite{RvW} for a two-phase flow. 

The starting point in this work is the Cahn--Hilliard--Navier--Stokes model developed in \cite{RvW}, which is describing the processes at the pore scale. The aim is to derive an upscaled model corresponding to the Darcy scale. We consider the simplified geometry of a thin strip, and assume that the  ratio of the width of the strip and its width is small. We employ asymptotic expansion methods that use this ratio as expansion parameter, and derive upscaled equations for transversally averaged quantities. In this respect, we follow the ideas in \cite{Carina15, KNP, Tycho1} for one-phase flow including precipitation and dissolution effects at the pore walls, and \cite{AndroStrip, AndroTube, Lunowa, Sohely} for two-phase flow, all considering a thin strip or tube. 
Observe that the pore-scale models in these works mentioned above are involving free boundaries. Instead, for the phase-field, pore scale model in \cite{CarinaLarsSorin} describing the flow of one fluid phase but including precipitation and dissolution, a Darcy-scale model is also derived for a thin strip by transversal averaging, in comparison to the one obtained by homogenization in more general situations.

This paper is organized as follows. First, in Section \ref{sec:si} a sharp-interface model for two fluid phases and one solid phase (including precipitation and dissolution) is presented. This model is approximated by the phase-field model proposed by \cite{RvW}, which is discussed briefly in Section \ref{sec:pf}. After bringing the phase-field model to a nondimensional form in Section \ref{sec:nondim}, in Section \ref{sec:up} we derive its upscaled counterpart by considering a thin strip geometry. The upscaled model still uses phase-field variables to locate the diffuse interfaces. In Section \ref{sec:sil} we identify the sharp-interface limit, that is the limit when letting the diffuse interface width go to zero. Notably the upscaling and the sharp-interface limit commute. The numerical examples discussed in Section \ref{sec:num} conclude the work.

\section{The Sharp-Interface Model}
\label{sec:si}
We start by presenting the sharp-interface model, which is then approximated by a phase-field model. We let $T > 0$ stand for the maximal time. For each $t \in [0, T]$, an $N$-dimensional domain $\Omega$ ($N = 2$ or 3) is partitioned into three disjoint subdomains, $\Omega_1(t)$, $\Omega_2(t)$ and $\Omega_3(t)$. These are occupied by the two fluid phases and by the solid phase respectively. The interface between the domain $\Omega_i$ and $\Omega_j$ is denoted by $\Gamma_{ij}$ ($i, j \in \set{1, 2, 3}, i \ne j$). Observe that these interfaces also depend on time.  

With $t \in (0, T]$, in the fluid occupied subdomains $\Omega_i(t)$, $i\in\set{1,2}$ the model is governed by the incompressible Navier--Stokes equations 
\begin{align*}
\nabla \cdot \vv &= 0, \\
\partial_t (\rho_i \vv) + \nabla \cdot \left( \rho_i \vv \otimes \vv \right) + \nabla p &= \nabla \cdot (2 \gamma_i \nabla^s \vv),
\end{align*}
where $\rho_i, \gamma_i$ denote the density, respectively viscosity of the  fluid phase $i$, all assumed constant here. $\vv$ and $p$ denote the fluid velocity and pressure in $\Omega_i$, the index $i$ being skipped. The symmetrized strain (Jacobian) is given by $\nabla^s \vv = \frac{1}{2}\left(\nabla \vv + (\nabla \vv)^t\right)$. 

At the interface $\Gamma_{12}(t)$ (separating $\Omega_1(t)$ and $\Omega_2(t)$) we assume that the velocity $\vv$ is continuous and that the jump in the normal stress is only in the normal direction, and proportional to the curvature of the interface
\begin{align*}
\jump{\vv} &= 0, \\
\jump{(p I - 2 \gamma \nabla^s \vv) \cdot \nv} &= \sigma_{12} \kappa \nv, \\
\nu &= \vv \cdot \nv. 
\end{align*}
Here $\jump{\cdot}$ denotes the jump of a quantity over the interface, $\nv$  the unit normal vector pointing outwards $\Omega_1$ and $\kappa$ the curvature of the interface. Through the last condition, the the normal velocity~$\nu$ of the interface and the normal velocity of the fluids are equal.

The subdomain $\Omega_3(t)$ is occupied by a mineral, formed by the precipitation of two solute species present in fluid 1. The reverse process, in which the mineral can be dissolved and release solute in fluid 1 is also possible. In a simplified setting, assuming a constant electrical charge, it suffices to consider only one solute concentration in the model, see \cite{VDuijnKnabner}, which is denoted by $c$. Here we assume that solute is only present in fluid 1. Therefore, the solute transport is governed by the transport-diffusion equation in $\Omega_1(t)$ 
\begin{align*}
\partial_t c + \nabla \cdot (\vv c) - D \Delta c = 0, 
\end{align*}
where $D$ is the constant diffusion coefficient. 

The interface $\Gamma_{13}(t)$ is evolving due to precipitation and dissolution.  At $\Gamma_{13}$ one has
\begin{align}
\nu = - r(c) + \alpha \sigma_{13} \kappa, \label{eq_SI_reaction}\\
D \nabla c \cdot \nv = \nu (c^\ast - c).  \label{eq_SI_reaction2}
\end{align}
The reaction rate $r(c)$ appearing in the former is generic. It accounts for dissolution and precipitation effects and is assumed increasing in $c$. The last term in \eqref{eq_SI_reaction}, involving a  constant $\alpha \geq 0$,  allows for curvature effects in the evolution of $\Gamma_{13}$. The latter is the Rankine-Hugoniot condition, ensuring the conservation of mass. Here  $c^\ast$ is a constant, similar to the concentration of the species as part of the mineral present in $\Omega_3$. Equations \eqref{eq_SI_reaction} and \eqref{eq_SI_reaction2} only hold at $\Gamma_{13}(t)$ and not at outer boundaries of $\Omega$. That is, we do not allow for precipitation and dissolution at the outer boundaries of $\Omega$.

At the fluid-fluid interface $\Gamma_{12}(t)$, a similar condition is imposed 
\begin{align*}
\nabla c \cdot \nv = 0.
\end{align*}
As before, $\nv$ is the unit normal vector pointing outwards $\Omega_1$. This follows from the Rankine-Hugoniot condition, since the concentration in fluid 2 is zero, and the normal velocity of the two fluids and of the interface are equal. 

In contrast to $\Gamma_{13}$, no precipitation or dissolution are possible at the interface $\Gamma_{23}$ between $\Omega_2$ and $\Omega_3$. This is because we assume that fluid 2 does not contain any solute species. Therefore, the interface does not evolve, and its normal velocity is $\nu = 0$.

Finally, at the interfaces between a fluid and the mineral a Navier-slip condition \cite{navier1823memoire} is assumed,  
\begin{align}
\label{eq_SI_slip}
\vv \cdot \tauv = -2 L_\text{slip} \tauv (\nabla^s \vv) \nv   
\end{align}
at  $\Gamma_{i3}$ ($i\in\set{1,2}$). Here $\tauv \in \RR^N$ is any tangent vector to $\Gamma_{i3}$ (thus $\tauv \perp \nv$). Here $L_\text{slip} \geq 0$ is a given slip length.

\section{The Phase-Field Model}
\label{sec:pf}
The sharp-interface model in Section \ref{sec:si} involves free boundaries, which makes it difficult from both analysis and numerical point of view. Relying on the idea to approximate the characteristic functions of each of the phases by smooth phase indicators \cite{Caginalp88}, phase-field models are convenient alternatives. For the specific problem considered here, a phase-field model called \textit{$\delta$-$2f1s$-model} was introduced in \cite{RvW}; here we present it briefly for completeness. We refer to  \cite{RvW} for more details on the derivation and the properties of the model, including the derivation of the sharp-interface limit. 

\subsection{Preliminaries}
The $\delta$-$2f1s$-model introduces three phase-field variables $\phi_1$, $\phi_2$, $\phi_3$ that represent the volume fraction of the two fluid phases and of the solid phase, respectively. Thus, $\phi_i$ approximates the indicator function of $\Omega_i$ appearing in the sharp-interface model in Section \ref{sec:si}. The phase-field variables $\Phiv = (\phi_1, \phi_2, \phi_3)^t$ are smooth and defined on the entire domain $\Omega$. In the sharp-interface model, the transition from one phase to another is across an interface. In the phase-field model, this interface is replaced by a diffuse transition zone from one phase to another, where the gradients of the corresponding phase-field variables are high. A ternary Cahn--Hilliard equation governs the evolution of $\Phiv$, and is coupled with a Navier--Stokes equation for fluid flow, and a reaction-transport-diffusion equation for dissolved ion concentration $c$. 

The $\delta$-$2f1s$-model additionally introduces a small regularisation parameter $\delta > 0$. Since no maximum principle holds for the Cahn--Hilliard equation, $\delta$ is used to ensure the positivity of the volume fractions. Also, the double-well potential
\begin{align}
\label{modelWdiph}
W_{\text{dw}}(\phi) = 450 \phi^4 (1-\phi)^4 + \delta \ell\left(\frac{\phi}{\delta} \right) + \delta \ell \left( \frac{1-\phi}{\delta} \right) , \text{ with } \ell(x) = \begin{cases} \frac{x^2}{1+x} & x \in (-1,0), \\ 0 & x \geq 0  \end{cases}
\end{align}
is employed. Observe that $W_\text{dw}$ has two minima at $0$ and $1$, and becomes unbounded at $-\delta$ and $1+\delta$. 
With this, we define the triple-well potential 
\begin{align}
\label{eq_W_W0}
W(\Phiv) := W_0(P \Phiv), \text{ where } W_0(\Phiv) = \sum_{i=1}^3 \Sigma_i W_\text{dw}(\phi_i). 
\end{align}
Here $\Sigma_i > 0$ are surface energy coefficients, and $P$ is the projection of $\RR^3$ onto the plane $\sum_i \phi_i = 1$, given by
\begin{align}
\label{eqProjection}
P\Phiv = \Phiv + \Sigma_T(1-\phi_1-\phi_2-\phi_3) \begin{pmatrix}
\Sigma_1^{-1}\\
\Sigma_2^{-1}\\
\Sigma_3^{-1}
\end{pmatrix}
,\qquad \frac{1}{\Sigma_T} = \frac{1}{\Sigma_1}+\frac{1}{\Sigma_2}+\frac{1}{\Sigma_3}.
\end{align}
As shown in \cite{RvW}, this construction ensures that the volume fractions sum up to one, i.e. $\sum_{i=1}^3 \phi_i = 1$, provided the initial data has this property. Furthermore, \cite{RvW} uses an energy argument and the unboundedness of the potential to show that  $-\delta < \phi_i < 1 + \delta$ ($i = 1, 2, 3$).

Next, we define the total fluid volume fraction $\tilde \phi_f$ and ion-dissolving fluid fraction $\phi_c$ as 
\begin{align}
&\tilde \phi_f := \phi_1+\phi_2 + 2\delta \phi_3, \label{rmodel_phif}\\
&\phi_c := \phi_1 \\
&\tilde \phi_c := \phi_1 + \delta, \label{rmodel_phic}
\end{align}
Here the tilde denotes a modification using the small parameter $\delta$, to ensure that the respective variables are positive. 
Using the (constant) fluid densities $\rho_i$ and viscosities $\gamma_i$ the total fluid density $\rho_f$ and viscosity $\tilde \gamma$ become 
\begin{align}
&\rho_f(\Phiv) := \rho_1 \phi_1 + \rho_2 \phi_2, \\ 
&\tilde \rho_f(\Phiv) := \rho_1 \phi_1 + \rho_2 \phi_2 + (\rho_1+\rho_2) \delta, \label{rmodel_rhof}\\ 
&\tilde \gamma(\Phiv) := \left(\phi_1 \gamma_1^{-1}+ \phi_2 \gamma_2^{-1} + \phi_3 \gamma_3^{-1}  + (\gamma_1^{-1}+\gamma_2^{-1}+\gamma_3^{-1})\delta\right)^{-1}, \label{rmodel_gamma}
\end{align}
As explained in \cite{RvW}, $\gamma_3$ is not the viscosity of the solid phase, but is chosen instead to archive a slip length $L_\text{slip}$ in the slip condition \eqref{eq_SI_slip}.

\subsection{The $\delta$-$2f1s$-Model}
We now present the $\delta$-$2f1s$-model. All equations are defined in $(0, T] \times \Omega$. The flow is governed by the Navier--Stokes equations and involves the fluid fraction $\tilde \phi_f$,
\begin{align}
\nabla \cdot (\tilde \phi_f \vv) &= 0, \label{rmodel1}\\
\begin{split}
 \partial_t (\tilde\rho_f \vv) + \nabla \cdot ((\rho_f \vv + \rho_1 \Jv_1 + \rho_2 \Jv_2) \otimes \vv) &= - \tilde\phi_f \nabla p + \nabla \cdot (2 \tilde \gamma(\Phiv) \nabla^s \vv) \\
 &\qquad - \rho_3 d(\tilde\phi_f) \vv + \tilde \Sv + \frac 12 \rho_1 \vv R.  \label{rmodel2}
 \end{split}
\end{align}
This is coupled with the transport-diffusion-reaction equation for the ion concentration 
\begin{align}
 \partial_t (\tilde\phi_c c) + \nabla \cdot ((\phi_c \vv + \Jv_1) c) &= D \nabla \cdot ( \tilde\phi_c  \nabla c) + c^\ast R. \label{rmodel3}
\end{align}
The phase-field variables $\phi_1$, $\phi_2$, $\phi_3$ are satisfying the Cahn--Hilliard equations 
\begin{align}
\partial_t \phi_1 + \nabla \cdot (\phi_1 \vv + \Jv_1) &= R, &&\hspace{-3em}   \label{rmodel4}\\
\partial_t \phi_2 + \nabla \cdot (\phi_2 \vv + \Jv_2) &= 0, &&\hspace{-3em}   \label{rmodel4a}\\
\partial_t \phi_3 + \nabla \cdot (2\delta \phi_3 \vv + \Jv_3) &= -R,  \label{rmodel4b}\\
\Jv_i &= - \frac{\eps M}{\Sigma_i} \nabla \mu_i, &&\hspace{-3em}  i\in\set{1,2,3},
\label{modelJ}\\
\mu_i &= \frac{\partial_{\phi_i}W(\Phiv)}{\eps} -  \eps\Sigma_i \Delta \phi_i, &&\hspace{-3em}  i\in\set{1,2,3}. \label{rmodel5}
\end{align}
Compared to the common Navier--Stokes equations, some modifications appear in \eqref{rmodel2}. The fluid density $\tilde \rho_f(\Phiv)$ introduces a strong coupling between the Navier--Stokes equations and the Cahn--Hilliard equations. All terms except the advection term use the modified quantities $\tilde \phi_f$, $\tilde \rho_f$ and $\tilde \gamma$. Additional flux terms $\rho_i \Jv_i \otimes \vv$ are introduced to account for momentum fluxes due to the Cahn--Hilliard evolution. Secondly, the dissipative term $-\rho_3 d(\tilde \phi_f) \vv$ is added. Here $d$ is a decreasing function s.t. $d(0) = d_0 > 0$ and $d(1)=0$, for example $d(\tilde \phi_f) = d_0 (1-\tilde \phi_f)^2$. The term $d(\tilde \phi_f)$ is therefore active in the solid phase and guarantees that $\vv$ remains small there. It also influences the slip length $L_\text{slip}$. Lastly, the surface tension term $\tilde \Sv$ is given by
\begin{align}
\label{rmodelS}
\tilde \Sv &= -\mu_2 \tilde \phi_f\nabla\left(\frac{\phi_1}{\tilde\phi_f}\right) -\mu_1 \tilde\phi_f\nabla\left(\frac{\phi_2}{\tilde\phi_f}\right)  - 2\delta \phi_3 \nabla (\mu_3-\mu_1-\mu_2) .
\end{align}

The reaction term $R$ modelling  precipitation and dissolution of ions is given by
\begin{align}
\label{modelrR}
\qquad R  = -q(\Phiv)  \left(r(c)+ \tilde \alpha \mu_1-\tilde \alpha \mu_3\right).
\end{align}
Here $r(c)$ is the increasing reaction rate used in the sharp interface description \eqref{eq_SI_reaction}. Additionally the precipitation process can depend on curvature effects through surface effects that are similar to surface diffusion, and are encountered if $\alpha > 0$. Again, the tilde denotes a modification of $\alpha$, that is $\tilde \alpha = \alpha + \delta$. Finally, to concentrate the reaction inside the diffuse interface region between fluid phase 1 and the solid phase, which is equivalent to the assumption made in the sharp-interface model, the non-dimensional term $q(\Phiv) = 30 \phi_1^2 \phi_3^2$ is used. Observe that $q$ dominates wherever neither $\phi_1$ nor $\phi_2$ are close to 0, which is precisely the envisaged location for the fluid 1 - mineral interface.

\section{Nondimensionalization}
\label{sec:nondim}
\begin{figure}
\centering
\includegraphics[width = 0.6\textwidth]{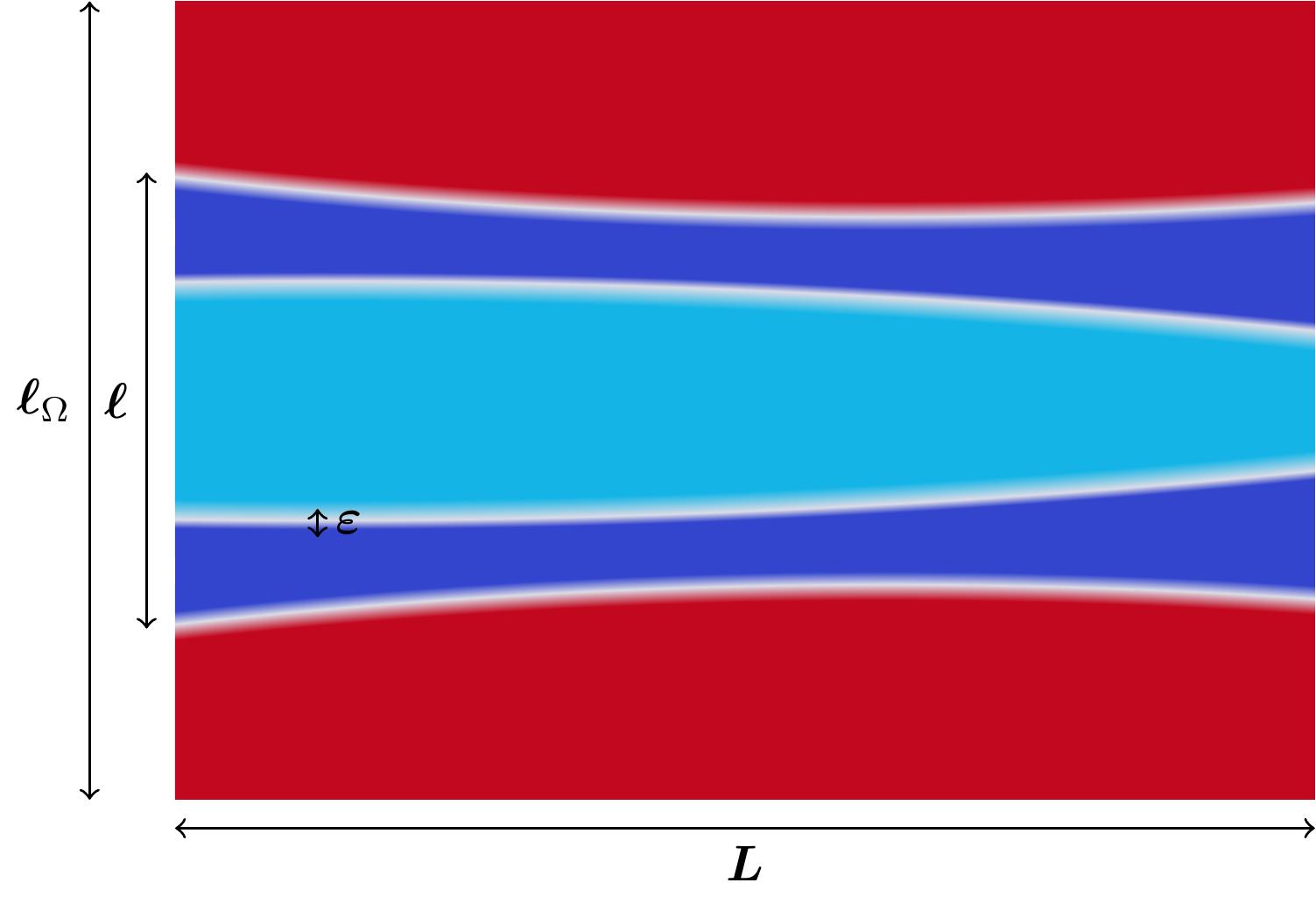}
\caption{Setting of the thin strip: The strip with length $L$ and width $\ell$ consists of solid walls (red, $\Phiv \approx (0,0,1)^t$) and fluid phases (light blue, dark blue). The diffuse interface width $\eps$ is smaller than $\ell$.}
\label{Figure_Scales}
\end{figure}

We proceed by bringing the $\delta$-$2f1s$-model \eqref{rmodel1}-\eqref{rmodel5} to a non-dimensional form, and derived an upscaled counterpart of it by employing asymptotic expansion and averaging techniques. We consider a simplified geometric setting. We start by introducing a thin strip having length $L$ and width $\ell \ll L$, as shown in figure~\ref{Figure_Scales}.

With a chosen domain width $\ell_\Omega > \ell$, the domain $\Omega = [0,L] \times [-\ell_\Omega/2, \ell_\Omega/2]$ includes the thin strip mentioned above, which is identified as $[0, L] \times [-\ell/ 2, \ell/ 2]$. The region outside the strip is occupied by the mineral, so $\Phiv \approx (0,0,1)^t$ there. The diffuse interfaces are located in regions of width $\eps$. We assume here that the diffuse-interface regions remain clearly separated inside the thin strip, hence $\eps \ll \ell$. 

Three length scales can be identified, $L\gg\ell\gg\eps$. These are related through the aspect ratio $\beta = \ell/L$ and the Cahn-Number $Cn = \eps/L$, both assumed small. Observe that, in fact, $Cn \ll \beta \ll 1$. 

The reference quantities used in the nondimensionalization procedure are listed in table \ref{table_nondim}. Nondimensional values are then identified by a hat. Note that we relate only few reference values directly to each other. In particular we do relate reference values when we want to emphasize an explicit dependence on $\r y$, as seen for $\r p$, $\r d$ and $\r \mu$. The choices are motivated as follows. To obtain an upscaled macroscopic velocity of order $\r v = \r x / \r t$, the pressure drop in the thin strip has to scale with $1/(\r y)^2$. Also, the slip length $L_\text{slip}$  is supposed to be of order $\ell$ and not $L$, which is achieved by a momentum dissipation scaling $1/(\r y)^2$. 
\begin{table}
\begin{tabular}{l|ll}
Variable                      & Reference value                        & Non-dimensional variable \\ \hline
time                           & $\r t = T$                                 & $\hat t = t/\r t$     \\
space                          & $\r x = L$,                            & $\hat x = x/\r x$     \\
                               & $\r y=\ell$,                           & $\hat y = y/\r y$     \\
                               & $\r \eps=\eps$                         & $\hat \eps = 1$       \\
velocity                       & $\r v = \r x / \r t$                   & $\hat \vv = \vv/\r v$ \\
density                        & $\r \rho = \rho_1$                     & $\hat \rho_i = \rho_i/\r \rho$, $i\in\set{1,2,3}$ \\
                               &                                        & $\hat \rho_f = \rho_f/\r \rho$ \\
                               &                                        & $\hat{\tilde \rho}_f = \tilde \rho_f/\r \rho$ \\
viscosity            & $\r \gamma = \gamma_1$                    & $\hat \gamma_i = \gamma_i/\r \gamma$, $i\in\set{1,2,3}$ \\
                               &                                        & $\hat {\tilde \gamma} = \tilde \gamma/\r \gamma$ \\
pressure                       & $\r p = \r \gamma \r v \r x/ (\r y)^2$ & $\hat p = p/\r p$     \\
momentum dissipation rate      & $\r d = \r \gamma / (\r \rho \r y^2)$
                        & $\hat d = d/\r d$     \\

surface energy                 & $\r \Sigma = \Sigma_1$                 & $\hat \Sigma_i = \Sigma_i/\r \Sigma$, $i\in\set{1,2,3}$\\
CH mobility        & $\r M = M$                             & $\hat M = 1$          \\
CH chemical potential       & $\r \mu = \r \Sigma / \r y$   & $\hat \mu = \mu/\r \mu$ \\
CH triple-well potential       & $\r W = \r \Sigma$         & $\hat W = W/\r \Sigma$\\
molar concentration            & $\r c = c^\ast$                        & $\hat c = c/\r c$     \\
diffusion coefficient          & $\r D = D$                             & $\hat D = 1$          \\
reaction rate                  & $\r r$                                 & $\hat r(\hat c) = r(c)/\r r$ \\
interface-reaction diffusivity\; & $\r \alpha = \r r / \r \mu$            & $\hat \alpha = \alpha /\r \alpha$\\                      
\end{tabular}
\caption{Variables, Reference Values and non-dimensional quantities for the nondimensionalization.}
\label{table_nondim}
\end{table}

We rewrite the Cahn number introduced above in terms of reference quantities, and define other dimensionless numbers that are used below to relate the reference quantities: the Reynolds number, Capillary number, Damk\"ohler number and P\'eclet numbers of the phase field and ion concentration, 
\begin{align}
\begin{alignedat}{3}
Re &= \frac{ \r \rho \r v \r x}{\r \gamma},& \qquad
Ca &= \frac{ \r \gamma \r v}{\r \Sigma}, & \qquad
Cn &= \frac{ \r \eps}{\r x},\\
Da &= \frac{ \r r \r x}{\r v}, & \qquad
Pe_{CH} &= \frac{ \r v \r x}{\r M}, & \qquad
Pe_{c} &= \frac{ \r v \r x}{\r D}. \label{eq_nondim_numbers}
\end{alignedat}
\end{align}
Clearly, the nondimensionalization also affects the spatial and temporal derivatives, namely  
\begin{align}
\hat \nabla = \r x \nabla, \quad \text{ and } \quad 
\partial_{\hat t} = \r t \partial_t. \label{eq_nondim_dt}
\end{align}

We now can insert the non-dimensional variables of table \ref{table_nondim}, the non-dimensional numbers \eqref{eq_nondim_numbers} and the non-dimensional operators in 
\eqref{eq_nondim_dt} into the $\delta$-$2f1s$-model \eqref{rmodel1}-\eqref{rmodel5}. The non-dimensional equations become  
\begin{align}
&\hat \nabla \cdot (\tilde \phi_f \hat \vv) = 0, \label{nmodel1}\\
\begin{split}
&\partial_{\hat t} (\hat {\tilde\rho}_f \hat \vv) + \hat \nabla \cdot (\hat \rho_f \hat \vv \otimes \hat \vv) + \frac{Cn}{\beta Pe_{CH}} \hat \nabla \cdot ((\hat \rho_1 \hat \Jv_1 + \hat \rho_2 \hat \Jv_2) \otimes \hat \vv) \\
 &\qquad = - \frac{1}{\beta^2 Re}\tilde\phi_f \hat\nabla \hat p +  \frac{1}{Re} \hat \nabla \cdot (2 \hat{\tilde \gamma}(\Phiv) \hat \nabla^s \hat \vv) \\
 &\qquad \qquad - \frac{1}{\beta^2 Re}\hat \rho_3 \hat d(\tilde\phi_f) \hat \vv + \frac{1}{\beta Re} \frac{1}{Ca} \hat {\tilde \Sv} +  Da \frac 12 \hat \rho_1 \hat \vv \hat R, \label{nmodel2}
 \end{split}
\end{align}
 for the flow, 
\begin{align}
 \partial_{\hat t} (\tilde\phi_c \hat c) + \hat \nabla \cdot (\phi_c \hat \vv \hat c) + \frac{Cn}{\beta Pe_{CH}} \hat \nabla \cdot(\hat \Jv_1 \hat c) &= \frac{1}{Pe_c} \hat \nabla \cdot ( \tilde\phi_c  \hat \nabla \hat c) + Da \hat R,\label{nmodel3}
 \end{align}
for the  ion transport-diffusion-reaction, while for the Cahn--Hilliard evolution one gets
\begin{align}\partial_{\hat t} \phi_1 + \hat\nabla \cdot (\phi_1 \hat \vv) + \frac{Cn}{\beta Pe_{CH}} \hat \nabla \cdot \hat \Jv_1 &= Da \hat R, &&\hspace{-3em}   \label{nmodel4}\\
\partial_{\hat t} \phi_2 + \hat\nabla \cdot (\phi_2 \hat \vv) + \frac{Cn}{\beta Pe_{CH}} \hat \nabla \cdot \hat \Jv_2 &= 0, &&\hspace{-3em}   \label{nmodel4a}\\
\partial_{\hat t} \phi_3 + \hat\nabla \cdot (2\delta \phi_3 \hat \vv) + \frac{Cn}{\beta Pe_{CH}} \hat \nabla \cdot \hat \Jv_3 &= -Da \hat R,  \label{nmodel4b}\\
\hat \Jv_i &= - \frac{1}{\hat \Sigma_i} \hat \nabla \hat \mu_i, &&i\in\set{1,2,3}, \label{nmodelJ}\\
\frac{\hat \mu_i}{\beta } &= \frac{\partial_{\phi_i}\hat W(\Phiv)}{Cn} -  Cn \hat \Sigma_i \hat \nabla^2 \phi_i, &&i\in\set{1,2,3}. \label{nmodel5}
\end{align}
All equations are defined in the dimensionless time-space domain $(0, 1] \times \hat{\Omega}$, where $\hat \Omega = [0,1] \times [- \hat \ell_\Omega/2, \hat \ell_\Omega/2]$. The 
surface tension and reaction are given as
\begin{align*}
\hat{\tilde \Sv} &= -\hat \mu_2 \tilde \phi_f\hat \nabla\left(\frac{\phi_1}{\tilde\phi_f}\right) -\hat \mu_1 \tilde\phi_f\hat \nabla\left(\frac{\phi_2}{\tilde\phi_f}\right)  - 2\delta \phi_3 \hat \nabla (\hat \mu_3-\hat \mu_1-\hat \mu_2) \\
\hat R &= -q(\Phiv) (\hat r(\hat c) + \hat {\tilde \alpha} \hat \mu_1 - \hat {\tilde \alpha} \hat \mu_3)
\end{align*}
From here on, we will only work with the non-dimensional model and therefore the hats are left out in the notation.

\section{Upscaling in a Thin Strip}
\label{sec:up}
We now proceed by deriving the upscaled model, obtained when passing to the limit $\beta \to 0$. This means that the thin strip reduces to a one-dimensional object, as its width is vanishing compared to its length.

We introduce new coordinates $(x,y)$ such that $\xv = (x, \beta y)$. In the thin strip we expect all variables to vary in longitudinal direction $\ev_x$ on the length scale $L = \r x$ and in transversal direction $\ev_y$ on the length scale $\ell = \r y = \beta \r x$. In particular this will result in $\nabla = \ev_x \partial_x + \beta^{-1} \ev_y \partial_y$.
\subsection{Scaling of Non-dimensional Numbers}
The upscaled model will also depend on the scaling of the dimensionless numbers \eqref{eq_nondim_numbers} with respect to $\beta$. We consider the following behavior of these  numbers with respect to $\beta$
\begin{align}
Re &= \overbar{Re} \label{eq_nr_re}\\
Ca &= \overbar{Ca} \label{eq_nr_ca}\\
Cn &= \beta \overbar \eps\label{eq_nr_cn}\\
Da &= \overbar{Da} / \overbar \eps \label{eq_nr_da}\\
Pe_{CH} &= 1 / (\beta^2  \overbar M) \label{eq_nr_pech}\\
Pe_{c} &= \overbar{Pe_c}\label{eq_nr_pec}
\end{align}
where $\overbar{Re}$, $\overbar{Ca}$, $\overbar \eps$, $\overbar{Da}$, $\overbar M$, $\overbar{Pe_c}$ are constants independent of $\beta$. In detail, these choices are motivated as follows. 
\begin{itemize}
\item The moderate Reynolds number \eqref{eq_nr_re} leads to a parabolic flow profile in the thin strip, we expect laminar flow. 
\item As the curvature of the fluid-fluid interface is of order $O(\beta)$, choosing a moderate capillary number $Ca$ in \eqref{eq_nr_ca} leads to the same pressure in both fluids, thus the capillary pressure becomes 0 (for sharp-interface models see also \cite{Lunowa, Sohely}). Note that this is a major difference to the three dimensional case, see e.g. \cite{AndroTube}, where we expect a curvature of $O(\beta^{-1})$ leading to a nonzero capillary pressure.
\item The scaling of the Cahn number $Cn$ in \eqref{eq_nr_cn} can be reformulated to $\overbar \eps = \eps / \ell$. Therefore the the interface width $\eps$ scales with the width of the thin strip, $\ell$. At the same time, the diffuse interface regions are assumed to be localised inside the thin strip, therefore we require $\eps \ll \ell$. This translates into a fixed, small $\overbar \eps$, i.e. $\overbar \eps \ll 1$. In the numerical experiments presented in Section \ref{sec:num} we choose $\overbar \eps = 0.03$.
\item We consider a moderate Damk\"ohler number \eqref{eq_nr_da}. 
In the sharp-interface model, this would ensure that the interfaces move with moderate velocity inside the thin strip, proportional to $\ell/T$. In the diffuse-interface model, the reaction is only active in the diffuse-interface region, which has an area scaling with $\eps$. Therefore $Da$ is divided by $\overbar \eps$, and expect to have fluid-solid or fluid-fluid interfaces evolving over the length scale $\ell$. A dominating Damk{\" o}hler regime like $Da = O(\beta^{-1})$ would instead lead to equilibrium-type reactions in the upscaled model, but the evolution of the interfaces should remain moderate. This can be achieved by assuming that the molar density of the species in the precipitate is sufficiently high to compensate the fast reaction kinetics. 
\item The high P\'eclet number \eqref{eq_nr_pech} for the phase field assures that the evolution of the phase field remains within the transversal length scale $\ell$ in an $O(1)$ timescale.
\item The moderate P\'eclet number of the ion diffusion \eqref{eq_nr_pec} will result in a macroscopic diffusion of ions, while the ion distribution in transversal direction equilibrates faster than the $O(1)$ timescale. 
\end{itemize}  
Lastly, the small, non-dimensional number $\delta>0$ appears in the $\delta$-$2f1s$-model. It is used as a regularisation parameter, to ensure the positivity of volume fractions, density and viscosity. Here we assume that $\delta$ is constant and independent of $\beta$.

\subsection{Asymptotic Expansions}
\label{sec:up_asy}
We assume that we can write solutions to the non-dimensional $\delta$-$2f1s$-model \eqref{nmodel1}-\eqref{nmodel5} in terms of an asymptotic expansion in $\beta$ of $\Phiv$, $\vv$, $p$, $c$, $\mu_1$, $\mu_2$, $\mu_3$. To be precise, we assume expansions of the form
\begin{align*}
\Phiv(t,\xv) &= \Phiv_0(t,x,y) + \beta \Phiv_1(t,x,y) + \beta^2 \Phiv_2(t,x,y) + \ldots\;,
\end{align*}
where $\Phiv_k$, $k \in \NN_0$ does not depend on $\beta$. In particular, we also use this notation for other variables, e.g.
\begin{align*}
\tilde \phi_{f} &= \phi_{f,0} + \beta \phi_{f,1} + \ldots
= \left(\phi_{1,0}+\phi_{2,0} + 2\delta \phi_{3,0} \right) + \beta \left( \phi_{1,1} + \phi_{2,1} + 2\delta \phi_{3,1}\right) + \ldots\;.
\end{align*}
Inserting these asymptotic expansions into the non-dimensional $\delta$-$2f1s$-model we group by powers of $\beta$. We use Taylor expansions to handle nonlinearities, e.g.
\begin{align*}
r(c) = r(c_0 + \beta c_1 + \ldots ) = r(c_0) + \beta r'(c_0) c_1 + O(\eps^2).
\end{align*}

\begin{remark}
\label{rem:scaling}
Note that the asymptotic expansions are written depending on the new coordinates $x$ and $y$. This means that in the $\ev_x$ direction variables can not vary on the (non-dimensional) length scale $\beta$, because a non-trivial function $f(x/\beta)$ can not be expanded in the form $f(x/\beta) = f_0(x) + \beta f_1(x) + \ldots$. In particular this implies that there are no phase-field interfaces possible perpendicular to the thin strip, as they would change the value of $\Phiv$ over the length $Cn = \beta \overbar \eps$. We will discuss in Section \ref{sec:num_shock} a numerical example that violates this assumption. 

The assumption is also violated for triple points, where all three phases meet, and for points where interfaces meet the boundary of $\Omega$ at $y = \pm \ell_\Omega/2$. Therefore $\ell_\Omega$ has to be chosen big enough, such that the width of the thin strip does not reach $\ell_\Omega$.
\end{remark}

The nondimensional domain is given by $\Omega = [0,1]\times [-\ell_\Omega/2, \ell_\Omega/2]$ and we choose as boundary conditions at $y = \pm \ell_\Omega/2$ for the upscaling, in detail
\begin{align}
\partial_y \Phiv(t,x,\pm \ell_\Omega/2) &= 0 \label{bc_phi}\\
\partial_y \mu(t,x,\pm \ell_\Omega/2) &= 0 \label{bc_mu}\\
\partial_y c(t,x,\pm \ell_\Omega/2) &= 0 \label{bc_c}\\
\vv(t,x,\pm \ell_\Omega/2) &= 0 \label{bc_v} 
\end{align}

\paragraph{Expansion of \eqref{nmodel1}, $O(\beta^{-1})$:}
Recall that $\nabla = \ev_x \partial_x + \beta^{-1} \ev_y \partial_y$. Therefore the leading order terms of \eqref{nmodel1} are of order $O(\beta^{-1})$, we have
\begin{align*}
\partial_y (\tilde \phi_{f,0} \vv_0) \cdot \ev_y = 0.
\end{align*}
We will denote components of $\vv$ as $\vv^{(1)} = \vv \cdot \ev_x$ and $\vv^{(2)} = \vv \cdot \ev_y$. Note that $\tilde \phi_{f,0}>0$ by construction in \eqref{rmodel_phif}, so after integrating and using the leading order of boundary condition \eqref{bc_v} we can divide by $\tilde \phi_{f,0}$ and obtain
\begin{align}
\label{v2=0}
\vv^{(2)}_0 = 0
\end{align}
As expected, there is no leading order flow perpendicular to the thin strip.

\paragraph{Expansion of \eqref{nmodel1}, $O(1)$:}
With \eqref{v2=0} we get in first order
\begin{align}
\label{eq_v0_v1}
\partial_x (\tilde \phi_{f,0} \vv^{(1)}_0) + \partial_y (\tilde \phi_{f,0} \vv^{(2)}_1)= 0.
\end{align}
The $O(\beta)$ term of boundary condition \eqref{bc_v} reads $\vv_1(y=\pm \ell_\Omega/2) = 0$. After integrating \eqref{eq_v0_v1} in $y$ we can use this to get
\begin{align}
\label{u_v0}
\partial_x \int_{-\ell_\Omega/2}^{\ell_\Omega/2} \tilde \phi_{f,0} \vv^{(1)}_0 \diff y = 0.
\end{align}
Here, $\tilde \phi_{f,0} \vv^{(1)}_0$ is the flux in $\ev_x$ direction, so \eqref{u_v0} implies that the total flux in $\ev_x$ direction is conserved.

\paragraph{Expansion of \eqref{nmodel5}, $O(\beta^{-1})$:}
We get with $Cn = \beta \overbar \eps$ three terms in leading order
\begin{align}
\label{u_mu}
\mu_{i,0} = \frac{\partial_{\phi_i} W(\Phiv_0)}{\overbar \eps} - \overbar \eps \Sigma_i \partial_y^2 \phi_{i,0}.
\end{align}
Notably from the Laplacian only derivatives in $\ev_y$-direction remain. In the upscaled model this will lead to a Cahn--Hilliard evolution that is only acting in $\ev_y$ direction.

\paragraph{Expansion of \eqref{nmodel4},\eqref{nmodel4a}, \eqref{nmodel4b}, $O(1)$:}
Note that with \eqref{eq_nr_cn}, \eqref{eq_nr_da} and \eqref{eq_nr_pech} we can write \begin{align}
\label{eq_nr_cnpech_dacn}
\frac{Cn}{\beta Pe_{CH}} = \beta^2 \overbar \eps \overbar M \qquad \text{and} \qquad Da = \frac{\overbar{Da}}{\overbar{\eps}}.
\end{align} 
We insert \eqref{nmodelJ} into \eqref{nmodel4},\eqref{nmodel4a}, \eqref{nmodel4b}, as we do not treat $\Jv_i$ as a primary variable. Together with \eqref{v2=0} we have in leading order $O(1)$
\begin{align}
\partial_{t} \phi_{1,0} + \partial_x (\phi_{1,0} \vv_0^ {(1)}) + \partial_y  (\phi_{1,0} \vv_1^ {(2)}) - \frac{\overbar \eps \overbar M}{\Sigma_1}  \partial_y^2 \mu_{1,0}  &= \frac{\overbar{Da}}{\overbar \eps } R_0, &&\hspace{-3em} \label{u_phi1}  \\
\partial_{t} \phi_{2,0} + \partial_x (\phi_{2,0} \vv_0^ {(1)}) + \partial_y  (\phi_{2,0} \vv_1^ {(2)}) - \frac{\overbar \eps \overbar M}{\Sigma_2}  \partial_y^2 \mu_{2,0}  &= 0, &&\hspace{-3em}  \label{u_phi2} \\
\partial_{t} \phi_{3,0} + \partial_x (2\delta \phi_{3,0} \vv_0^ {(1)}) + \partial_y  (2\delta \phi_{3,0} \vv_1^ {(2)}) - \frac{\overbar \eps \overbar M}{\Sigma_3}  \partial_y^2 \mu_{3,0}  &= - \frac{\overbar{Da}}{\overbar \eps }R_0, \label{u_phi3}
\end{align}
where the leading order term of the reaction is given by
\begin{align}
R_0 = - q(\Phiv_0)(r(c_0) + \tilde \alpha \mu_{1,0} - \tilde \alpha \mu_{3,0}).\label{R0}
\end{align}
Note that as in \eqref{u_mu} only the $y$-derivatives of the Laplacian remain in the leading order.

\paragraph{Expansion of \eqref{nmodel3}, $O(\beta^{-2})$:}
We obtain in leading order only one $O(\beta^{-2})$ term
\begin{align*}
\frac{1}{\overbar{Pe_c}} \partial_y ( \tilde \phi_{c,0} \partial_y c_0) = 0
\end{align*}
Integrating in $y$ and using the leading order term of boundary condition \eqref{bc_c} results in
\begin{align*}
\tilde \phi_{c,0} \partial_y c_0 = 0
\end{align*}
Because by construction $\tilde \phi_{c,0} > 0$, we conclude 
\begin{align}
\label{u_c0}
\partial_y c_0=0.
\end{align}
Therefore $c_0$ is constant in $\ev_y$ direction, and we write $c_0 = c_0(t,x)$ to emphasize that $c_0$ only depends on the $x$ coordinate.

\paragraph{Expansion of \eqref{nmodel3}, $O(\beta^{-1})$:}
As we found $\partial_y c_0=0$ in \eqref{u_c0}, we get in first order only the term
\begin{align*}
\frac{1}{\overbar{Pe_c}} \partial_y ( \tilde \phi_{c,0} \partial_y c_1) = 0
\end{align*}
With analogous argumentation to the $O(\beta^{-2})$ case we get $\partial_y c_1=0$ and can write $c_1 = c_1(t,x)$ to show that $c_1$ is independent of $y$.

\paragraph{Expansion of \eqref{nmodel3}, $O(1)$:}
Similar to the $O(1)$ expansion of \eqref{nmodel4},\eqref{nmodel4a},\eqref{nmodel4b}, we insert the Cahn--Hilliard flux $\Jv_i$ \eqref{nmodelJ} and the non-dimensional numbers \eqref{eq_nr_cnpech_dacn} into the equation, and use \eqref{v2=0}. We obtain the second order terms
\begin{align*}
&\partial_{t} (\tilde\phi_{c,0} c_0) + \partial_x (\phi_{c,0} \vv^{(1)}_0 c_0) + \partial_y (\phi_{c,0} \vv^{(2)}_1 c_0) - \frac{\overbar \eps \overbar M}{\Sigma_1} \partial_y( c_0 \partial_y \mu_{1,0}) \\
&\quad = \frac{1}{\overbar{Pe_c}} \partial_x ( \tilde \phi_{c,0} \partial_x c_0) + \frac{1}{Pe_c} \partial_y ( \tilde \phi_{c,0} \partial_y c_2) + \frac{\overbar{Da}}{\overbar \eps } R_0
\end{align*}
where $R_0$ is given by \eqref{R0}. After integrating in $y$ we can use the boundary conditions \eqref{bc_mu}, \eqref{bc_c}, \eqref{bc_v} to eliminate the terms containing a $y$ derivative. We obtain
\begin{align}
\label{u_c}
&\frac{d}{dt} \left( c_0  \int_{-\ell_\Omega/2}^{\ell_\Omega/2} \tilde\phi_{c,0} \diff y \right) + \partial_x \left(c_0  \int_{-\ell_\Omega/2}^{\ell_\Omega/2} \phi_{c,0} \vv^{(1)}_0 \diff y \right) \\
&\qquad= \frac{1}{\overbar{Pe_c}} \partial_x \left( \left( \int_{-\ell_\Omega/2}^{\ell_\Omega/2} \tilde \phi_{c,0} \diff y \right) \partial_x c_0\right) + \frac{\overbar{Da}}{\overbar \eps }  \int_{-\ell_\Omega/2}^{\ell_\Omega/2} R_0 \diff y
\end{align}
Here we have written $c_0$ outside of the integrals to emphasize that $c_0$ does not depend on $y$. Equation \eqref{u_c} is a transport-diffusion-reaction equation for $c_0(x,t)$, where the coefficients still depend on the exact distribution of $\Phiv_0$  in the $\ev_y$ direction.

\paragraph{Expansion of \eqref{nmodel2}, $O(\beta^{-3})$:}
The only term of order $O(\beta^{-3})$ is
\begin{align*}
-\frac{1}{\overbar{Re}} \tilde \phi_{f,0} \ev_y \partial_y p_0= 0
\end{align*}
As $\tilde \phi_{f,0}$ is positive by construction, we conclude that $p_0$ does not depend on $y$ and write $p_0 = p_0 (t,x)$.

\paragraph{Expansion of \eqref{nmodel2}$\cdot \ev_x$, $O(\beta^{-2})$:}
We investigate in the first order only the equation for the $x$-component. With \eqref{v2=0} and $p_0 = p_0(t,x)$ the remaining terms are
\begin{align*}
- \frac{1}{\overbar{Re}} \tilde\phi_{f,0} \partial_x p_0 + \frac{1}{\overbar{Re}} \partial_y (\tilde \gamma(\Phiv_0) \partial_y \vv^{(1)}_0) - \frac{1}{\overbar{Re}} \rho_3 d(\tilde\phi_{f,0}) \vv^{(1)}_0 = 0
\end{align*}
We can interpret this a a linear differential equation for $\vv^{(1)}_0$ with boundary conditions $\eqref{bc_v}$. In particular we can use the linearity to write
\begin{align}
\label{eq_w_v}
\vv^{(1)}_0(t,x,y) = -w(t,x,y) \partial_x p_0(t,x)
\end{align}
where $w$ is the solution to the cell problem
\begin{align}
\label{eq_w}
\rho_3 d(\tilde\phi_{f,0}) w - \partial_y (\tilde \gamma(\Phiv_0) \partial_y w) &= \tilde\phi_{f,0},\\
\label{eq_w_bc}
w(t,x, \pm \ell_\Omega/2) &= 0
\end{align}
For a given $\Phiv$ the function $w$ calculates the parabolic flow profile in the cross section of the thin strip. As we expect from a Darcy-type flow, the fluid velocity is proportional to $-\partial_x p_0$, shown in \eqref{eq_w_v}.

\begin{remark}
We note that by construction $\tilde \gamma > 0$ and therefore the cell problem \eqref{eq_w}, \eqref{eq_w_bc} has a unique solution.
\end{remark}

\subsection{Upscaling in a Thin Strip: Summary}
Let us summarize the results of the upscaling. Except for $\vv$ we will only need the leading order term of each unknown, and will therefore drop the subscript $0$. We will call the model \eqref{umodel1}-\eqref{umodel_v2} the upscaled $\delta$-$2f1s$-model.

From \eqref{u_v0} and \eqref{eq_w_v} we have the macroscopic continuity equation for the total flux $Q_f$ and the Darcy-equation for the pressure $p$, and the macroscopic transport-diffusion-reaction equation for the ion concentration $c$ \eqref{u_c}
\begin{align}
\label{umodel1}
\partial_x Q_f &= 0, \\
\label{umodel1a}
Q_f &= -K_f \partial_x p\\
\label{umodel3}
\frac{d}{dt} \left( \tilde \phi_{c,\text{total}} c \right) + \partial_x \left( (-K_c \partial_x p) c \right) &= \frac{1}{\overbar{Pe_c}} \partial_x \left( \tilde \phi_{c,\text{total}} \partial_x c\right) + \frac{\overbar{Da}}{\overbar \eps } R_\text{total}
\end{align}
These equations are macroscopic in the sense that the unknowns $Q_f$, $p$ and $c$ depend only on $x$ and $t$, but not on $y$. The parameters in these equations are upscaled quantities, depending on the exact distribution of the phases in $y$ direction 
\begin{align}
\label{umodel_phictot}
\tilde \phi_{c,\text{total}} &=   \int_{-\ell_\Omega/2}^{\ell_\Omega/2} \tilde \phi_{c} \diff y \\
\label{umodel_kf}
K_f(t,x) &=  \int_{-\ell_\Omega/2}^{\ell_\Omega/2} \tilde \phi_{f} w \diff y \\
K_c(t,x) &=  \int_{-\ell_\Omega/2}^{\ell_\Omega/2} \tilde \phi_{c} w \diff y \\
\label{umodel_rtot}
R_\text{total} &=  \int_{-\ell_\Omega/2}^{\ell_\Omega/2} R \diff y
\end{align}

For the phase-field parameters we still have to solve the fully coupled 2-d problem \eqref{u_mu}, \eqref{u_phi1}, \eqref{u_phi2}, \eqref{u_phi3}, that is
\begin{align}
\partial_{t} \phi_{1} + \partial_x (\phi_{1} \vv_0^ {(1)}) + \partial_y  (\phi_{1} \vv_1^ {(2)}) - \frac{\overbar \eps \overbar M}{\Sigma_1}  \partial_y^2 \mu_{1}  &= \frac{\overbar{Da}}{\overbar \eps }R, &&\hspace{-3em} \label{umodel4}  \\
\partial_{t} \phi_{2} + \partial_x (\phi_{2} \vv_0^ {(1)}) + \partial_y  (\phi_{2} \vv_1^ {(2)}) - \frac{\overbar \eps \overbar M}{\Sigma_2}  \partial_y^2 \mu_{2}  &= 0, &&\hspace{-3em}  \label{umodel4a} \\
\partial_{t} \phi_{3} + \partial_x (2\delta \phi_{3} \vv_0^ {(1)}) + \partial_y  (2\delta \phi_{3} \vv_1^ {(2)}) - \frac{\overbar \eps \overbar M}{\Sigma_3}  \partial_y^2 \mu_{3}  &= -\frac{\overbar{Da}}{\overbar \eps }R, \label{umodel4b} \\
\mu_{i} = \frac{\partial_{\phi_i} W(\Phiv)}{\overbar \eps} - \overbar \eps \Sigma_i \partial_y^2 \phi_{i} \qquad  i\in\set{1,2,3} \label{umodel5}
\end{align}
with the reaction term
\begin{align}
\label{umodel_r}
R = - q(\Phiv)(r(c) + \tilde \alpha \mu_{1} - \tilde \alpha \mu_{3})
\end{align}
Note that in contrast to the non-dimensional model \eqref{nmodel1}-\eqref{nmodel5} the Cahn--Hilliard evolution acts only in $\ev_y$ direction. The only term acting in $\ev_x$ direction is the transport of the fluid phases. This will enable us in Section \ref{sec:num_up} to develop a numerical model that uses explicit upwinding for the fluid transport and can therefore decouple cell-problems for different values of $x$.

For the flow it suffices to solve the cell problem \eqref{eq_w}, \eqref{eq_w_bc}
\begin{align}
\label{umodel2}
\rho_3 d(\tilde\phi_{f}) w - \partial_y (\tilde \gamma(\Phiv) \partial_y w) &= \tilde\phi_{f},\\
\label{umodel2_bc}
\lim_{y\to\pm {\ell_\Omega/2}} w &= 0
\end{align}
and recover the flow $\vv^{(1)}_0$, $\vv^{(2)}_1$ by \eqref{eq_w_v} and \eqref{eq_v0_v1}
\begin{align}
\label{umodel_v1}
\vv^{(1)}_0 &= -w \partial_x p \\
\label{umodel_v2}
\partial_x (\tilde \phi_{f} \vv^{(1)}_0) + \partial_y (\tilde \phi_{f} \vv^{(2)}_1) &= 0
\end{align}

\section{Sharp-Interface Limit of the Upscaled $\delta$-$2f1s$-Model}
\label{sec:sil}
In the previous section we have investigated the scale separation $\beta = \ell / L \to 0$. A different limit process that is commonly investigated for phase-field models is the sharp-interface limit $\eps \to 0$. In \cite{RvW} this limit is analyzed for the $\delta$-$2f1s$ model \eqref{rmodel1}-\eqref{rmodel5}, resulting in the sharp-interface evolution described in Section \ref{sec:si}.

Because the upscaled $\delta$-$2f1s$-model \eqref{umodel1}-\eqref{umodel_v2} still contains a Cahn--Hilliard evolution, depending on the small number $\overbar \eps = \eps / \ell$, we can investigate the sharp-interface limit $\overbar \eps \to 0$ of the upscaled $\delta$-$2f1s$-model. This means that we are interested in the limit process of vanishing diffuse interface width $\eps$ compared to the width $\ell$ of the thin strip. In the following we will use matched asymptotic expansions to analyze this limit, the argumentation is mostly analogous to \cite{RvW}.

\subsection{Assumptions and Scaling of Non-dimensional Numbers}
To derive the the sharp-interface limit $\overbar \eps \to 0$, we assume that $\overbar{Pe_c}, \overbar{Da}, \overbar{M}$ are constant and independent of $\overbar \eps$. This choice of scaling allows for a reasonable limit process, with physical properties independent of the diffuse interface width.

The scaling $\delta = \overbar \eps$ is important. The regularisation parameter $\delta$ is introduced in the $\delta$-$2f1s$ to ensure the positivity of e.g. the density $\tilde \rho_f(\Phiv)$ in \eqref{rmodel_rhof}. This $\delta$-regularisation is not necessary for the sharp-interface formulation, and the choice $\delta = \overbar \eps$ leads to $\delta$ vanishing in the sharp-interface limit. 

As a basic assumption we expect to have solutions that form bulk phases, characterized by nearly constant $\Phiv$, and interfaces, characterized by a large gradient of $\Phiv$. We also assume that $\mu_i$, $i\in\set{1,2,3}$ is of order $O(1)$, not of order $O(\overbar \eps^{-1})$, as equation \eqref{umodel5} would suggest. For a discussion of why this assumption is reasonable on a $O(1)$ timescale, see \cite{pego}.

We also assume that in an interface between phase $\Phiv = \ev_i$ and $\Phiv = \ev_j$ the third phase is not present. This assumption is reasonable because with our constructions of $W$ \eqref{eq_W_W0} minimizers of the Ginzburg-Landau energy $W(\Phiv) + \sum_i \Sigma_i \Delta \phi_i$ that connect $\Phiv = \ev_i$ and $\Phiv = \ev_j$ satisfy $\phi_k = 0$, $k\in\set{1,2,3}\setminus\set{i,j}$.

\subsection{Outer Expansions}
\label{sec:sil_outer}
For the bulk phases we assume that we can write solutions to the upscaled $\delta$-$2f1s$-model \eqref{umodel1}-\eqref{umodel_v2} in terms of an outer asymptotic expansion in $\overbar \eps$ for the variables $\Phiv$, $w$, $\vv^{(1)}_0$, $\vv^{(2)}_1$, $p$, $c$, $\mu_1$, $\mu_2$, $\mu_3$. That is, similar to the expansions in Section \ref{sec:up_asy}, we assume expansions of the form
\begin{align*}
\Phiv(t,x,y) = \Phiv^{out}_0(t,x,y) + \overbar \eps \Phiv^{out}_1(t,x,y) + \overbar \eps^2 \Phiv^{out}_2(t,x,y) + \ldots
\end{align*} 
Here the outer expansion terms $\Phiv^{out}_k$, $k\in \NN_0$ are independent of $\overbar \eps$. The expansions for the macroscopic variables $p(x), c(x)$ do not depend on $y$. We will insert these expansions into the upscaled $\delta$-$2f1s$-model and group by orders of $\bar \eps$. Analogous to Section \ref{sec:up} we handle nonlinearities by Taylor expansion.

\paragraph{Outer Expansion of \eqref{umodel5}, $O(\beta^{-1})$:}
We can argue as in \cite{RvW} to find that the only stable solutions to the leading order terms are $\Phiv^{out}_0 = \ev_k$, $k\in\set{1,2,3}$ with the restriction $\phi_{k,1}^{out} \leq 0$ and $\phi_{i,1}^{out}, \phi_{j,1}^{out}\geq 0$ for $\set{i,j}=\set{1,2,3}\setminus \set{k}$. The additional restriction stems from the fact that the triple well potential $W$ depends on $\delta = \overbar \eps$.

We define the set $\Omega_k(t)$ to be the set of $(x,y)$ where $\Phiv^{out}_0(t,x,y) = \ev_k$. In the sharp interface formulation $\Omega_k(t)$ will represent the domain of phase $k$.

\paragraph{Outer Expansion of \eqref{umodel2}, $O(1)$:}
In $\Omega_3$, i.e. in case $\Phiv^{out}_0 = \ev_3$, we have $\tilde \phi^{out}_{f,0} = 0$ and the leading order reads
\begin{align}
\label{eq_wout}
\rho_3 d_0 w^{out}_0 - \partial_y (\gamma_3 \partial_y w^{out}_0) = 0
\end{align}
where $d_0 = d(0)>0$. In the fluid phases $\Omega_i$, $i\in\set{1,2}$, we have $\Phiv^{out}_0 = \ev_i$ and therefore $\tilde \phi^{out}_{f,0} = 1$. Note that by construction $d(1)=0$. With this we obtain in leading order
\begin{align}
\label{eq_woutb}
 - \partial_y (\gamma_i \partial_y w^{out}_0) = 1
\end{align}

\paragraph{Outer Expansion of \eqref{umodel_v2}, $O(1)$:}
In the fluid phases $\Phiv^{out}_0 = \ev_i$, $i\in\set{1,2}$ we have  $\tilde \phi^{out}_{f,0} = 1$ and obtain
\begin{align}
\label{eq_dxv_out}
\partial_x (\vv^{(1),out}_{0,0}) + \partial_y (\vv^{(2),out}_{1,0}) &= 0
\end{align}

\paragraph{Outer Expansion of \eqref{umodel1}, \eqref{umodel1a}, $O(1)$:}
We now consider the macroscopic equations. The equations for the flow \eqref{umodel1}, \eqref{umodel1a} upscale trivially, the leading order reads
\begin{align}
\label{eq_dxQ}
\partial_x Q^{out}_{f,0} &= 0, \\
\label{eq_dxP}
Q^{out}_{f,0} &= -K^{out}_{f,0} \partial_x p^{out}_0
\end{align}
where the parameter $K^{out}_{f,0}$ is the leading order expansion of $K_f$, using \eqref{umodel_kf}
\begin{align}
\label{eq_Koutf0}
K^{out}_{f,0} =  \int_{-\ell_\Omega/2}^{\ell_\Omega/2} \phi^{out}_{f,0} w^{out}_{f,0} \diff y.
\end{align}
Note that the leading order expansion of $\tilde \phi_f$ is $\phi^{out}_{f,0}$ as the $\delta$-modification is of order $O(\overbar \eps)$ because of the scaling choice $\delta = \overbar \eps$.

\paragraph{Outer Expansion of \eqref{umodel3}, $O(1)$:}
For the transport-diffusion-reaction equation for $c$ let us first investigate the reaction term. We have with \eqref{umodel_rtot} and \eqref{umodel_r}
\begin{align*}
\frac{\overbar{Da}}{\overbar \eps } R_\text{total} = - \frac{\overbar{Da}}{\overbar \eps }  \int_{-\ell_\Omega/2}^{\ell_\Omega/2} q(\Phiv)(r(c) + \tilde \alpha \mu_{1} - \tilde \alpha \mu_{3}) \diff y
\end{align*}
As $q(\Phiv^{out}) = O(\eps^2)$ in the bulk phases $\phi^{out}_0 = \ev_k$, $k\in\set{1,2,3}$, there is no contribution of the reaction term in the bulk at leading order. Note that there will be a contribution of this term in the interface regions, see Section \ref{sec:sil_inner}. Overall we have for \eqref{umodel3} in leading order
\begin{align}
\label{eq_cout}
\frac{d}{dt} \left( \phi^{out}_{c,\text{total},0} c^{out}_0 \right) + \partial_x \left( (-K^{out}_{c,0} \partial_x p^{out}_0) c^{out}_0 \right) &= \frac{1}{\overbar{Pe_c}} \partial_x \left( \phi^{out}_{c,\text{total},0} \partial_x c^{out}_0\right) + \overbar{Da} R_{\text{interface},0}
\end{align}
with coefficients
\begin{align}
\label{eq_phioutc}
\phi^{out}_{c,\text{total},0} &=  \int_{-\ell_\Omega/2}^{\ell_\Omega/2} \phi^{out}_{c,0} \diff y,\\
\label{eq_Koutc0}
K^{out}_{c,0} &=  \int_{-\ell_\Omega/2}^{\ell_\Omega/2} \phi^{out}_{c,0} w^{out}_{f,0} \diff y
\end{align}
and $R_{\text{interface},0}$ as a placeholder for the interface contributions of the reaction term.
\subsection{Inner Expansions}
\label{sec:sil_inner}
We have shown in Section \ref{sec:sil_outer} that the domain is partitioned into $\Omega_1$, $\Omega_2$ and $\Omega_3$. We locate the interfaces between the phases as 
\begin{align}
\label{eq_gammaij}
\Gamma_{ij}(t) = \set{(x,y)\in \Omega : \phi_i(t,x,y) = \phi_j(t,x,y) \geq 1/3}.
\end{align}
We assume that $\Gamma_{ij}$ is a smooth, one-dimensional manifold. As explained in Remark \ref{rem:scaling} we do not consider triple-points, where all three phases meet, and do not allow for the interfaces to touch the boundary of $\Omega$ at $y=\pm\ell_\Omega/2$. Also, interfaces can not occur perpendicular to the thin strip and therefore there exists locally around an interface $\Gamma_{ij}$ a unique mapping $s(t,x)$ such that $(x,s(t,x)) \in \Gamma_{ij}$.

We use this mapping to introduce a new coordinate $z$ close to the interface
\begin{align*}
z(x,t) = \frac{y - s(t,x)}{\overbar \eps}.
\end{align*}
Because we expect the interface width to be of size $\overbar \eps$ , the coordinate $z$ is scaled by $\eps^{-1}$. The velocity of $\Gamma_{ij}$ at $(x,s)$ in $y$-direction is given by
\begin{align*}
\nu(x) = \partial_t s(t,x)
\end{align*}
We will use the new coordinates $(t,x,z)$ as the coordinates to describe the interfaces $\Gamma_{ij}$. For a generic function $f(t,x,y)=f^{in}(t,x,z)$ we obtain the transformation rules
\begin{align}
\label{eq_trafo_dt}
\partial_t f &= - \frac{1}{\eps} \nu \partial_z f^{in} + \partial_t f^{in} \\
\label{eq_trafo_dy}
\partial_y f &= \frac{1}{\eps} \partial_z f^{in} \\
\label{eq_trafo_dx}
\partial_x f &= - \frac{1}{\eps} (\partial_x s) \partial_z f^{in} + \partial_x f^{in}
\end{align} 

We assume that close to an interface $\Gamma_{ij}$ we can write solutions to the upscaled $\delta$-$2f1s$-model \eqref{umodel1}-\eqref{umodel_v2} in terms of an inner asymptotic expansion in $\overbar \eps$ for the variables $\Phiv$, $w$, $\vv^{(1)}_0$, $\vv^{(2)}_1$, $\mu_1$, $\mu_2$, $\mu_3$. That is we assume expansions of the form
\begin{align*}
\Phiv(t,x,y) = \Phiv^{in}_0(t,x,z) + \overbar \eps \Phiv^{in}_1(t,x,z) + \overbar \eps^2 \Phiv^{in}_2(t,x,z) + \ldots
\end{align*}
with coefficients $\Phiv^{in}_k$ independent of $\overbar \eps$. In contrast to the outer asymptotic expansions, the inner asymptotic expansions depending on the $(t,x,z)$ coordinates. This will lead to different terms being of highest order when inserting the expansions into the upscaled$\delta$-$2f1s$ model. We do not use inner expansions of the macroscopic variables $p$ and $c$, as they are constant across all interfaces.

To relate inner and outer expansions, we match the limit value of inner expansions for $z \to \pm \infty$ with the limit value of the outer expansions at $s$ (from the respective side). The matching conditions are well studied \cite{Caginalp88}, we use
\begin{align}
\label{matching0}
\lim_{z \to \pm \infty} \Phiv^{in}_0(t,x,z) &= \lim_{y\to 0^+} \Phiv^{out}_0(t,x,s\pm y) \\
\label{matchingD0}
\lim_{z \to \pm \infty} \partial_z \Phiv^{in}_0(t,x,z) &= 0 \\
\label{matchingD1}
\lim_{z \to \pm \infty} \partial_z \Phiv^{in}_1(t,x,z) &= \lim_{y\to 0^+} \partial_y \Phiv^{out}_0(t,x,s\pm y)
\end{align}

\paragraph{Inner Expansion of \eqref{umodel5}, $O(\overbar \eps^{-1})$:}
Consider an interface between bulk phases $\Phiv^{out}_0 = \ev_i$ and $\Phiv^{out}_0 = \ev_j$. With matching condition \eqref{matching0} this means
\begin{align}
\label{eq_phiv_pm}
\lim_{z\to -\infty}\Phiv^{in}_0 = \ev_i \qquad \text{and}\qquad \lim_{z\to \infty}\Phiv^{in}_0 = \ev_j.
\end{align}
Then by assumption we have no third phase contributions across the interface, that is
\begin{align*}
\phi^{in}_{k,0} = 0, \quad \text{with } k\in\set{1,2,3}\setminus\set{i,j}
\end{align*}
Following the argument in \cite{RvW} we calculate the leading order terms of \eqref{umodel5} for $\mu_k$ and find $\phi^{in}_{j,0}$ as a solution to the ordinary differential equation
\begin{align}
\label{eq_w-dzz}
W'_{\text{dw}} (\phi^{in}_{j,0}) - \partial_z^2 \phi^{in}_{j,0} = 0
\end{align}
with additional conditions 
\begin{align*}
\lim_{z\to -\infty}\phi^{in}_{j,0} = 0, \qquad \lim_{z\to -\infty}\phi^{in}_{j,0} = 1, \qquad \phi^{in}_{j,0}(t,x,0) = 1/2 .
\end{align*}
The first two conditions are boundary conditions from \eqref{eq_phiv_pm} while the third condition stems from definition of $\Gamma_{ij}$ \eqref{eq_gammaij} and centers the interface at $z=0$. With a lengthy calculation $\phi^{in}_{j,0}$ is implicitly given by
\begin{align}
\label{eq_phi_in}
z = \frac{1}{30}\left(\frac{1}{1-\phi^{in}_{j,0}}-\frac{1}{\phi^{in}_{j,0}} + 2 \log\left(\frac{\phi^{in}_{j,0}}{1-\phi^{in}_{j,0}} \right) \right).
\end{align}
We find $\phi^{in}_{i,0}$ by $\phi^{in}_{i,0} = 1- \phi^{in}_{j,0}$.

\paragraph{Inner Expansion of \eqref{umodel_v2}, $O(\overbar \eps^{-1})$:}
Using the coordinate transformations \eqref{eq_trafo_dy} and \eqref{eq_trafo_dx}, we get in leading order
\begin{align*}
-(\partial_x s) \partial_z (\phi^{in}_{f,0} \vv^{(1),in}_{0,0}) + \partial_z (\phi^{in}_{f,0} \vv^{(2),in}_{1,0}) = 0.
\end{align*}
Note that $\partial_x s(t,x)$ does not depend on $z$ and therefore 
\begin{align}
\label{eq_dxsv_const}
-(\partial_x s)\phi^{in}_{f,0} \vv^{(1),in}_{0,0} + \phi^{in}_{f,0}\vv^{(2),in}_{1,0} = const.
\end{align}
with respect to $z$. 
Across the interface $\Gamma_{12}$ we have $\phi^{in}_{f,0} = 1$ and thus with matching condition \eqref{matching0} we get for all $z\in\RR$ 
\begin{align}
\begin{split}
&-(\partial_x s) \vv^{(1),in}_{0,0}(t,x,z) + \vv^{(2),in}_{1,0}(t,x,z) \\
&\quad= \lim_{z\to \pm \infty} -(\partial_x s)\vv^{(1),in}_{0,0}(t,x,z) + \vv^{(2),in}_{1,0}(t,x,z) \\
&\quad= \lim_{y\to 0^+} -(\partial_x s)\vv^{(1),out}_{0,0}(t,x,s\pm y) + \vv^{(2),out}_{1,0}(t,x,s\pm y)
\label{eq_gamma12_dxsv}
\end{split}
\end{align}
In particular this means that the term $-(\partial_x s)\vv^{(1),out}_{0,0} + \vv^{(2),out}_{1,0}$ is continuous across the $\Gamma_{12}$ interface.

When matching \eqref{eq_dxsv_const} at the fluid-solid interfaces $\Gamma_{13}$ and $\Gamma_{23}$, $\phi^{in}_{f,0}$ vanishes in the limit towards the solid phase, we can conclude
\begin{align}
\label{eq_dxsv_0}
-(\partial_x s)\phi^{in}_{f,0}\vv^{(1),in}_{0,0} + \phi^{in}_{f,0}\vv^{(2),in}_{1,0} = 0
\end{align}
Using matching condition \eqref{matching0} we find
\begin{align}
\label{eq_gamma13_dxsv}
-(\partial_x s)\vv^{(1),out}_{0,0} + \vv^{(2),out}_{1,0} = 0
\end{align}
for the fluid velocity. This condition therefore allows only for fluid flow parallel to the fluid-solid interfaces.

\paragraph{Inner Expansion of \eqref{umodel4},\eqref{umodel4a},\eqref{umodel4b}, $O(\overbar \eps^{-1})$:}
We will argue analogous to \cite{RvW}. The leading order expansions for \eqref{umodel4}, \eqref{umodel4a} and \eqref{umodel4b} are given by
\begin{align}
\begin{split}
\label{eq_umodel4_leading}
&-\nu\partial_z \phi^{in}_{1,0} -(\partial_x s) \partial_z (\phi^{in}_{1,0} \vv^{(1),in}_{0,0}) + \partial_z (\phi^{in}_{1,0} \vv^{(2),in}_{1,0}) - \frac{\overbar{M}}{\Sigma_1}\partial_z^2 \mu^{in}_{1,0} 
\\
&\qquad= -\overbar{Da}\, q(\Phiv^{in}_0)(r(c^{out}_0) + \alpha \mu^{in}_{1,0} - \alpha \mu^{in}_{3,0})
\end{split}\\
\begin{split}
\label{eq_umodel4a_leading}
&-\nu\partial_z \phi^{in}_{2,0} -(\partial_x s) \partial_z (\phi^{in}_{2,0} \vv^{(1),in}_{0,0}) + \partial_z (\phi^{in}_{2,0} \vv^{(2),in}_{1,0}) - \frac{\overbar{M}}{\Sigma_2}\partial_z^2 \mu^{in}_{2,0} = 0
\end{split}\\
\begin{split}
\label{eq_umodel4b_leading}
&-\nu\partial_z \phi^{in}_{3,0} - \frac{\overbar{M}}{\Sigma_3}\partial_z^2 \mu^{in}_{3,0} = +\overbar{Da}\, q(\Phiv^{in}_0)(r(c^{out}_0) + \alpha \mu^{in}_{1,0} - \alpha \mu^{in}_{3,0})
\end{split}
\end{align}

Let us first consider the interface $\Gamma_{13}$, with $\Omega_1$ being in the negative $z$ direction. Here $\phi^{in}_{1,0} = \phi^{in}_{f,0}$ and with \eqref{eq_dxsv_0} the advection terms vanish from \eqref{eq_umodel4_leading}. We also have no third phase contributions and therefore $\phi^{in}_{1,0} + \phi^{in}_{3,0} = 1$. With notation $\mu_{3-1} := \mu^{in}_{3,0} - \mu^{in}_{1,0}$ we calculate $\Sigma_3\, \cdot$ \eqref{eq_umodel4b_leading} $- \Sigma_1\, \cdot$ \eqref{eq_umodel4_leading}
\begin{align}
\label{eq_mu3-1}
-(\Sigma_1+\Sigma_3) \nu \partial_z \phi^{in}_{3,0} - \overbar{M} \partial_z^2\mu_{3-1} = (\Sigma_1+\Sigma_3) \overbar{Da} \, q(\Phiv^{in}_0)(r(c^{out}_0) - \alpha \mu_{3-1})
\end{align}
In \cite{RvW} it is shown that with \eqref{eq_w-dzz} and by construction of $q$ the identity $q(\Phiv^{in}_0) = \partial_z \phi^{in}_{3,0}$ holds. We can interpret \eqref{eq_mu3-1} as an ordinary differential equation for $\mu_{3-1}$ with boundary conditions $\lim_{z\to\pm\infty} \partial_z \mu_{3-1} = 0$ (by using matching condition \eqref{matchingD0}). 

In the case $\alpha = 0$ all constant functions $\mu_{3-1}$ are solutions to the differential equation, under the compatibility condition 
\begin{align}
\label{eq_mu3-1_alpha0}
\nu = - \overbar{Da}\,  r(c^{out}_0).
\end{align}
In case $\alpha > 0$ the unique solution to \eqref{eq_mu3-1} is given by the constant function
\begin{align}
\label{eq_mu3-1_alpha}
\mu_{3-1} = \alpha^{-1} (\nu + \overbar{Da}\,  r(c^{out}_0)).
\end{align}
We can combine \eqref{eq_mu3-1_alpha0} and \eqref{eq_mu3-1_alpha}, and also consider the case that the fluid and solid side of the $\Gamma_{13}$ interface is switched. Overall we conclude
\begin{align}
\label{eq_nu13_alpha}
\nu = \begin{cases}
\alpha (\mu^{in}_{1,0}-\mu^{in}_{3,0})+\overbar{Da}\,r(c) & \text{ if } \lim_{z\to -\infty}\Phiv^{in}_0 = \ev_3 \text{ and }\lim_{z\to \infty}\Phiv^{in}_0 = \ev_1 \\
\alpha (\mu^{in}_{3,0}-\mu^{in}_{1,0})-\overbar{Da}\,r(c) & \text{ if } \lim_{z\to -\infty}\Phiv^{in}_0 = \ev_1 \text{ and }\lim_{z\to \infty}\Phiv^{in}_0 = \ev_3 
\end{cases}
\end{align}

For $\Gamma_{23}$ we can argue analogous to the $\Gamma_{13}$ case. Because there is no precipitation, i.e. $q(\Phiv^{in}_0) = 0$, we obtain
\begin{align}
\label{eq_nu23}
\mu^{in}_{3,0}-\mu^{in}_{2,0}= const. \qquad\text{and}\qquad \nu = 0 .
\end{align}

Lastly, we consider the fluid-fluid interface $\Gamma_{12}$, with $\Omega_1$ in the direction of negative $z$. There is no precipitation process, so with $q(\Phiv^{in}_0) = 0$ we integrate over \eqref{eq_umodel4_leading} and use matching conditions \eqref{matching0} for $\phi^{in}_{1,0}$ and \eqref{matchingD0} for $\partial_z \mu^{in}_{1,0}$ and obtain
\begin{align}
\label{eq_nu12}
\nu = - (\partial_x s) \vv^{(1),in}_{0,0} + \vv^{(2),in}_{1,0}
\end{align}
Furthermore $\mu^{in}_{1,0}$ has to be constant in $z$, and with analogous argumentation using \eqref{eq_umodel4a_leading} also $\mu^{in}_{2,0}$ is constant.

\paragraph{Inner Expansion of \eqref{umodel5}, $O(1)$:}
We consider the interface $\Gamma_{ij}$ with $\Omega_i$ in negative $z$ direction. We assume the absence of a third phase, that is $\phi^{in}_{k,0} = 0$, $k\in\set{1,2,3}\setminus\set{i,j}$, and find by construction of $W$ in \eqref{eq_W_W0} that $\partial_{\phi_i} W'(\Phiv^{in}_0) = W'_\text{dw}(\phi^{in}_{i,0})$.  We examine the difference $\mu_i - \mu_j$ at first order and find
\begin{align}
\label{eq_mui-j}
\mu^{in}_{i,0}-\mu^{in}_{j,0} = \Sigma_i W_{\text{dw}}''(\phi^{in}_{i,0})\phi^{in}_{i,1} - \Sigma_i \partial_z^2 \phi^{in}_{i,1} - \Sigma_j W_{\text{dw}}''(\phi^{in}_{j,0})\phi^{in}_{j,1} + \Sigma_j \partial_z^2 \phi^{in}_{j,1}
\end{align}
In absence of a third phase $\phi^{in}_{i,0}+\phi^{in}_{j,0} = 1$, and by construction $W_{\text{dw}}(\phi)$ is symmetric around $\phi=1/2$. Therefore $W_{\text{dw}}''(\phi^{in}_{i,0}) = W_{\text{dw}}''(\phi^{in}_{j,0})$, and we rewrite \eqref{eq_mui-j} as
\begin{align*}
\mu^{in}_{i,0}-\mu^{in}_{j,0} = \left(W_{\text{dw}}''(\phi^{in}_{j,0}) - \partial_z^2 \right) \left(\Sigma_i \phi^{in}_{i,1} - \Sigma_j \phi^{in}_{j,1} \right) \\
\end{align*}
Recall that $\mu^{in}_{i,0}-\mu^{in}_{j,0}$ is constant across the interface $\Gamma_{ij}$. After multiplying with $\partial_z \phi^{in}_{j,0}$ and integrating over $z$ we calculate
\begin{align*}
\mu^{in}_{i,0}-\mu^{in}_{j,0} &= \int_{-\infty}^{\infty} \left(\partial_z \phi^{in}_{j,0}\right)  \left(\mu^{in}_{i,0}-\mu^{in}_{j,0}\right) \diff z \\
&= \int_{-\infty}^{\infty} \left(\partial_z \phi^{in}_{j,0}\right) \left(W_{\text{dw}}''(\phi^{in}_{j,0}) - \partial_z^2 \right) \left(\Sigma_i \phi^{in}_{i,1} - \Sigma_j \phi^{in}_{j,1} \right) \diff z \\
&= \int_{-\infty}^{\infty} \left(W_{\text{dw}}''(\phi^{in}_{j,0}) \partial_z \phi^{in}_{j,0} - \partial_z^3 \phi^{in}_{j,0} \right) \left(\Sigma_i \phi^{in}_{i,1} - \Sigma_j \phi^{in}_{j,1} \right) \diff z \\
&= \int_{-\infty}^{\infty} \partial_z\left(W_{\text{dw}}'(\phi^{in}_{j,0})  - \partial_z^2 \phi^{in}_{j,0} \right) \left(\Sigma_i \phi^{in}_{i,1} - \Sigma_j \phi^{in}_{j,1} \right) \diff z \\
&= 0
\end{align*}
We have used partial integration to get to the third line, the boundary terms vanish with matching condition \eqref{matchingD1} and the structure of $\phi^{in}_{j,0}$ \eqref{eq_phi_in}. The fourth line evaluates to zero with the identity \eqref{eq_w-dzz}. Note that compared to \cite{RvW} there is no curvature term in this calculation, as the Cahn--Hilliard evolution acts only in the $y$-direction.

We conclude
\begin{align}
\label{eq_mui_muj}
\mu^{in}_{i,0}=\mu^{in}_{j,0}
\end{align}
and \eqref{eq_nu13_alpha} simplifies to
\begin{align}
\label{eq_nu13}
\nu = \begin{cases}
+\overbar{Da}\;r(c) & \text{ if } \lim_{z\to -\infty}\Phiv^{in}_0 = \ev_3 \text{ and }\lim_{z\to \infty}\Phiv^{in}_0 = \ev_1 \\
-\overbar{Da}\;r(c) & \text{ if } \lim_{z\to -\infty}\Phiv^{in}_0 = \ev_1 \text{ and }\lim_{z\to \infty}\Phiv^{in}_0 = \ev_3 
\end{cases}
\end{align}

\paragraph{Inner Expansion of \eqref{umodel2}, $O(\overbar \eps^{-2})$:}
At leading order the equation reads
\begin{align*}
\partial_z ( \gamma( \Phiv^{in}_0 ) \partial_z w^{in}_0) = 0
\end{align*}
After integrating in $y$ we use matching condition \eqref{matchingD0} divide by $\gamma( \Phiv^{in}_0 )>0$ and find
\begin{align}
\label{eq_dzw}
\partial_z w^{in}_0 = 0.
\end{align}
that is $w$ is constant across the interface. With matching condition \eqref{matchingD0} this implies
\begin{align}
\label{eq_wout_interface}
\lim_{y\to 0^+} w^{out}_0(t,x,s+y) = \lim_{y\to 0^+} w^{out}_0(t,x,s-y) .
\end{align}

\paragraph{Inner Expansion of \eqref{umodel2}, $O(\overbar \eps^{-1})$:}
With \eqref{eq_dzw} the first order term of \eqref{umodel2} reads
\begin{align*}
\partial_z ( \gamma( \Phiv^{in}_0 ) \partial_z w^{in}_1) = 0
\end{align*}
We integrate and with matching conditions \eqref{matching0}, \eqref{matchingD1} we get
\begin{align}
\begin{split}
\label{eq_woutgamma}
&\lim_{y\to 0^+} \left(\gamma(\Phiv^{out}_0(t,x,s+y) \partial_y w^{out}_0(t,x,s+y)\right) \\
&\qquad= 
\lim_{y\to 0^+} \left(\gamma(\Phiv^{out}_0(t,x,s-y) \partial_y w^{out}_0(t,x,s-y)\right) .
\end{split}
\end{align}

\paragraph{Inner Expansion of \eqref{umodel3}, $O(1)$:}
We only need to investigate the reaction term
\begin{align*}
\frac{\overbar{Da}}{\overbar \eps } R_\text{total} = - \frac{\overbar{Da}}{\overbar \eps }  \int_{-\ell_\Omega/2}^{\ell_\Omega/2} q(\Phiv)(r(c) + \tilde \alpha \mu_{1} - \tilde \alpha \mu_{3}) \diff y
\end{align*}
On $\Gamma_{12}$ and $\Gamma_{23}$ we have $q(\Phiv^{in}) = O(\overbar \eps^2)$ and therefore no leading order contribution. Let us consider $\Gamma_{13}$ with $\Omega_1$ in negative $z$ direction. Using \eqref{eq_mui_muj} the leading order term of the integrand is $q(\Phiv^{in}_0) r(c^{out}_0)$. Transforming the integral to the $z$ coordinate results in the leading order term of $O(1)$
\begin{align*}
- \overbar{Da}\; r(c^{out}_0) \int_{-\infty}^{\infty} q(\Phiv^{in}_0) \diff z
\end{align*}
In \cite{RvW} it is shown that by construction of $q$ we have $q(\Phiv^{in}_0) = \diff_z \phi^{in}_{3,0}$. With matching condition \eqref{matching0} the integral evaluates to one. When considering $\Gamma_{13}$ with $\Omega_1$ in positive $z$ direction we get the same result.

There might be multiple $\Gamma_{13}$ interfaces contributing to the macroscopic reaction term. Therefore, the total contribution to \eqref{umodel3} at order $O(1)$ is
\begin{align}
\label{eq_rinterface}
 \overbar{Da} \;R_\text{interface} &= - \overbar{Da} N(\Gamma_{13}) r(c^{out}_0)
\end{align}
with $N(\Gamma_{13})$ being the number of $\Gamma_{13}$ interfaces for a fixed $x$.

\subsection{Sharp-Interface Limit: Summary}
\label{sec:us}
We will summarize the results of the matched asymptotic expansions. We drop the subscript $0$ and the superscript $out$ for ease of notation. We call \eqref{usmodel1}-\eqref{usmodel_vi} the upscaled sharp-interface model.

The macroscopic equations for the unknowns $Q_f$, $p$ and $c$ are given by \eqref{eq_dxQ}, \eqref{eq_dxP} and \eqref{eq_cout}, that is
\begin{align}
\label{usmodel1}
\partial_x Q_f &= 0, \\
\label{usmodel1a}
Q_f &= -K_f \partial_x p\\
\label{usmodel3}
\frac{d}{dt} \left( \phi_{c,\text{total}} c \right) + \partial_x \left( (-K_c \partial_x p) c \right) &= \frac{1}{\overbar{Pe_c}} \partial_x \left( \phi_{c,\text{total}} \partial_x c\right) + \overbar{Da}\; R_\text{interface}
\end{align}

The coefficients of the upscaled equations depend on the distribution of the phases in the thin strip. In contrast to the upscaled phase-field model \eqref{umodel1}-\eqref{umodel_v2} the sharp-interface limit does not depend on the phase-field variables $\Phiv$. Instead the three disjoint domains $\Omega_1(t)$, $\Omega_2(t)$ and $\Omega_3(t)$ are used to locate the phases. The interface between $\Omega_i$ and $\Omega_j$ is denoted by $\Gamma_{ij}$. We introduce the notation $\Omega_i|_x = \set{y\in [-\ell_\Omega/2, \ell_\Omega/2] : (x,y)\in \Omega_i(t)}$, and write $N(\Gamma_{13})$ for the number of $\Gamma_{13}$ interfaces at a given $x$. With \eqref{eq_phioutc}, \eqref{eq_Koutf0}, \eqref{eq_Koutc0}, \eqref{eq_rinterface} we can calculate the coefficients of \eqref{usmodel1}-\eqref{usmodel3} as
\begin{align}
\phi_{c,\text{total}}(x) &=  \text{vol}\left(\Omega_1|_x\right)\\
 K_f(t,x) &= \int_{\Omega_1|_x \cup \Omega_2|_x} w \diff y \\
 \label{usmodel_kc}
 K_c(t,x) &= \int_{\Omega_1|_x} w \diff y \\
 R_\text{interface} &= - N(\Gamma_{13}) r(c)
\end{align}

We describe the evolution of the phases with the interface velocity $\nu$. This velocity in $y$ direction is given by \eqref{eq_nu13}, \eqref{eq_nu23}, \eqref{eq_nu12}, summarized as
\begin{align}
\label{usmodel_i}
\nu &= \pm \overbar{Da}\;r(c) && \text{ on } \Gamma_{13}, \text{ with } \Omega_1 \text{ in } \pm y \text{ direction} \\
\label{usmodel_i2a}
\nu &= 0 && \text{ on } \Gamma_{23} \\
\label{usmodel_i3}
\nu &= -(\partial_x s) \vv^{(1)}_0 + \vv^{(2)}_1 && \text{ on } \Gamma_{12}.
\end{align}

For the flow profile we solve at each $x$ and $t$ a cell problem for the unknown $w$. Summarizing \eqref{eq_wout},\eqref{eq_woutb}, \eqref{eq_wout_interface}, \eqref{eq_woutgamma} and the boundary condition \eqref{umodel2_bc}, the unknown $w$ is given by the second order differential equation
\begin{align}
\label{usmodel2}
 - \partial_y (\gamma_1 \partial_y w) &= 1 && \text{ in } \Omega_1|_x \\
 \label{usmodel2b}
 - \partial_y (\gamma_2 \partial_y w) &= 1 && \text{ in } \Omega_2|_x \\
 \label{usmodel2c}
\rho_3 d_0 w - \partial_y (\gamma_3 \partial_y w) &= 0 && \text{ in } \Omega_3|_x \\
\label{usmodel2_i}
\jump{w} &= 0 && \text{ at } \Gamma_{12}, \Gamma_{13}, \Gamma_{23}\\
\label{usmodel2_ib}
\jump{\gamma \partial_y w} &= 0 && \text{ at } \Gamma_{12}, \Gamma_{13}, \Gamma_{23} \\
\label{usmodel2_bc}
w &= 0 && \text{ at } y = \pm \ell_\Omega/2
\end{align}

For the transport of the fluid-fluid interface $\Gamma_{12}$ in \eqref{usmodel_i3} we need the flow velocities $\vv^{(1)}_0$ and $\vv^{(2)}_1$. We then get the horizontal flow velocity $\vv^{(1)}_0$ from \eqref{umodel_v1}, that is
\begin{align}
\label{usmodel_dxp}
\vv^{(1)}_0 = - w \partial_x p_0
\end{align}
For the vertical flow velocity $\vv^{(2)}_1$ one has to solve \eqref{eq_dxv_out},\eqref{eq_gamma12_dxsv} and \eqref{eq_gamma13_dxsv}, summarized
\begin{align}
\label{usmodel_v}
\partial_x (\vv^{(1)}_{0}) + \partial_y (\vv^{(2)}_{1}) &= 0 && \text{ in } \Omega_1 \cup \Gamma_{12} \cup \Omega_2 \\
\label{usmodel_vi}
-(\partial_x s) \vv^{(1)}_0 + \vv^{(2)}_{1} &= 0 && \text{ on } \Gamma_{13} \text{ and } \Gamma_{23}
\end{align}

\subsection{Upscaled Sharp-Interface Model in a Simplified Geometry with Symmetry}
\label{sec:us_geo}
The upscaled sharp-interface model \eqref{usmodel1}-\eqref{usmodel_vi} uses no assumption on how the phases are distributed. When these are appearing in a fixed order, the model simplifies. In this case, there is no need to consider a general subdomain $\Omega_i$ for the  phase $i$, it is sufficient to know the width of the phase $i$ layer in the $y$ direction. These widths become unknowns of the model.

\begin{figure}
\centering
\begin{tikzpicture}
\begin{scope}[scale = 3.0]
\draw[name path=ly0, dashed, -stealth] (0,-1) -- (0,1);
\node[right] at (0,1) {$y$};
\path[name path=ly1] (2,-1) -- (2,1);
\path[name path=ly2] (3,-1) -- (3,1);
\draw[name path=lx1] (0,-0.7) to[out=20,in=210] (2,-0.7) to[out=30,in=200] 
		(4,-0.7);
\draw[name path=lx2] (0,-0.4) to[out=20,in=195] 
		(4,-0.3);
\draw[name path=lx3, dashed, -stealth] (0,0) -- (4,0);
\node[right] at (4,0) {$x$};
\draw[name path=lx4] (0,0.4) to[out=-20,in=-195] 
		(4,0.3);
\draw[name path=lx5] (0,0.7) to[out=-20,in=-210] (2,0.7) to[out=-30,in=-200] 
		(4,0.7);
\path [name intersections={of=lx3 and ly1,by=E1}];
\path [name intersections={of=lx4 and ly1,by=E2}];
\path [name intersections={of=lx5 and ly1,by=E3}];
\path [name intersections={of=lx1 and ly2,by=F1}];
\path [name intersections={of=lx2 and ly2,by=F2}];
\path [name intersections={of=lx4 and ly2,by=F4}];
\path [name intersections={of=lx5 and ly2,by=F5}];
\node[inner sep=0pt] (d0) at (2,0) {};
\node[align=center] at (1,-0.85) {$\Omega_3$};
\node[fill=white,align=center] at (1,-0.5) { $\Omega_1$};
\node[fill=white,align=center] at (1,0) {$\Omega_2$};
\node[fill=white,align=center] at (1,0.5) { $\Omega_1$};
\node[align=center] at (1,0.85) {$\Omega_3$};
\node[fill=white,align=center] at (F1) {$\Gamma_{13}$};
\node[fill=white,align=center] at (F2) {$\Gamma_{12}$};
\node[fill=white,align=center] at (F4) {$\Gamma_{12}$};
\node[fill=white,align=center] at (F5) {$\Gamma_{13}$};
\draw[stealth-stealth] (E1) -- node[right]{$d_2(x,t)$} (E2);
\draw[stealth-stealth] (E2) -- node[right]{$d_1(x,t)$} (E3);
\end{scope}
\end{tikzpicture}  

 \caption{Symmetric geomery of two fluid phases in a thin strip}
  \label{Figure_Geometry}
\end{figure}
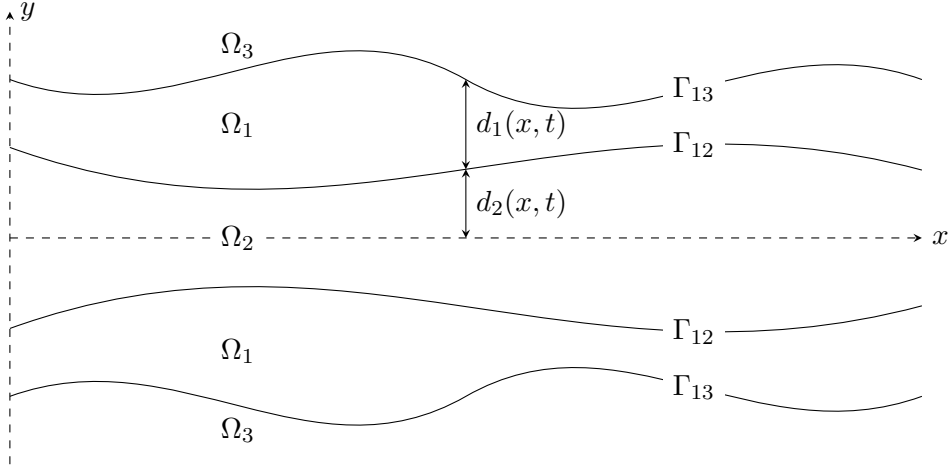

We assume here the following simplified geometry. The solid phase (in $\Omega_3$) is covered by a film of fluid $1$ (occupying $\Omega_1$). The second fluid (in $\Omega_2$) is located in the middle of the thin strip. For simplicity we assume symmetry around the $x$-axis. An illustration of the geometry is given in figure \ref{Figure_Geometry}.

With functions $d_1(t,x) > 0$, $d_2(t,x) > 0$, representing the width in $y$ direction of the fluid phase 1, respectively 2, we can describe this situation by defining
\begin{align*}
\Omega_2(t) =& \set{(x,y) : -d_2(t,x) < y < d_2(t,x) } \\
\Omega_1(t) =&\set{(x,y) : -d_1(t,x)-d_2(t,x) < y < -d_2(t,x)} \\
		   &\quad \cup \set{(x,y) : d_2(t,x) < y < d_1(t,x)+d_2(t,x) } \\
\Omega_3(t) =&\set{(x,y) : -\ell_\Omega/2 < y < -d_1(t,x)-d_2(t,x)} \\
		   &\quad \cup \set{(x,y) : d_1(t,x)+d_2(t,x)< y < \ell_\Omega/2 }
\end{align*}

In this geometry the solution $w$ to the cell problem \eqref{usmodel2}-\eqref{usmodel2_bc} depends only on the variables $d_1$ and $d_2$, and on the choice of $\ell_\Omega$. With a lengthy calculation we find that the terms depending on $\ell_\Omega$ decay exponentially fast for big $\ell_\Omega$, and we drop them in the following. The remaining terms lead to
\begin{align*}
K_f &= \frac{2}{\gamma_1}\left(\frac{(d_1+d_2)^3}{3} + \left(\frac{\gamma_1}{\gamma_2} - 1\right)\frac{d_2^3}{3} + L_\text{slip}(d_1+d_2)^2\right) \\
K_c &= \frac{2}{\gamma_1}\left(\frac{d_1^3}{3} + \frac{d_1^2 d_2}{2} + L_\text{slip} d_1(d_1+d_2)\right) \\
\end{align*}
with the slip length $L_\text{slip}$ given by
\begin{align*}
L_\text{slip} = \frac{\gamma_1}{\sqrt{\rho_3 d_0 \gamma_3}}.
\end{align*}

We can relate $\partial_t d_1$ and $\partial_t d_2$ with the interface velocities \eqref{usmodel_i}-\eqref{usmodel_i3}. Considering the fluid-solid interface $\Gamma_{13}$ we get with \eqref{usmodel_i}
\begin{align}
\label{eq_dtd1d2}
\partial_t \left( d_1 +  d_2 \right) &= \nu = - \overbar{Da} \; r(c)
\end{align}
while for the fluid-fluid interface $\Gamma_{12}$ we calculate with \eqref{usmodel_i3}, \eqref{usmodel_v}, \eqref{usmodel_vi}
\begin{align*}
\partial_t d_2 &= \nu \\
&= - (\partial_x d_2) \vv^{(1)}_0(t,x,d_2) + \vv^{(2)}_1(t,x,d_2) \\
&= - (\partial_x d_2) \vv^{(1)}_0(t,x,d_2) + \vv^{(2)}_1(t,x,d_2) \\
&\qquad + (\partial_x (d_1 + d_2) ) \vv^{(1)}_0(t,x,d_1+d_2) - \vv^{(2)}_1(t,x,d_1+d_2)\\
&= (\partial_x (d_2+d_1)) \vv^{(1)}_0(t,x,d_1+d_2) - (\partial_x d_2) \vv^{(1)}_0(t,x,d_2) -\int_{d_2}^{d_2+d_1} \partial_y \vv^{(2)}_1(t,x,y) \diff y \\
&= (\partial_x (d_2+d_1)) \vv^{(1)}_0(t,x,d_1+d_2) - (\partial_x d_2) \vv^{(1)}_0(t,x,d_2) +\int_{d_2}^{d_2+d_1} \partial_x \vv^{(1)}_0(t,x,y) \diff y \\
&= \partial_x \left( \int_{d_2}^{d_2+d_1} \vv^{(1)}_0(t,x,y) \diff y \right) 
\end{align*}
The integral equals the total fluid flux in $x$ direction in the upper half of $\Omega_1$. We use \eqref{usmodel_dxp}, \eqref{usmodel_kc} and the symmetry of $w$ around $y=0$ to further calculate
\begin{align}
\label{eq_dtd2}
\partial_t d_2 &= \partial_x \left( \int_{d_2}^{d_2+d_1} \vv^{(1)}_0 \diff y \right) 
 = -\partial_x \left( (\partial_x p) \int_{d_2}^{d_2+d_1} w \diff y  \right) = -\frac{1}{2} \partial_x\left(K_c \partial_x p\right)
\end{align}

We can now summarize \eqref{usmodel1},\eqref{usmodel1a},\eqref{usmodel3},\eqref{eq_dtd1d2} and \eqref{eq_dtd1d2} as an upscaled model for the unknowns $d_1$, $d_2$, $p$, $Q_f$ and $c$
\begin{align}
\label{xmodel1}
\partial_t d_1 + \partial_t d_2 &= - \overbar{Da} \; r(c(t,x)) \\
\label{xmodel2}
\partial_t d_2 &= - \frac{1}{2} \partial_x \left(K_c(d_1,d_2) \partial_x p\right)  \\
\label{xmodel3}
Q_f &= -K_f(d_1,d_2) \partial_x p \\
\label{xmodel4}
\partial_x Q_f &= 0 \\
\label{xmodel5}
\frac{d}{dt} \left( 2 d_1 c \right) + \partial_x \left( (-K_c(d_1,d_2) \partial_x p) c \right) &= \frac{1}{Pe_c} \partial_x \left( 2 d_1 \partial_x c\right) - 2 \overbar{Da}\; r(c) 
\end{align}

\begin{remark}
\label{rem:hyp}
We can rewrite \eqref{xmodel2},\eqref{xmodel3} to highlight the hyperbolicity of the model. As discussed in Remark \ref{rem:scaling} one assumption for the upscaling is that there is no occurrence of triple points. Therefore we assume $d_1 > 0$ and $d_2 > 0$ and deduce $K_f > 0$, $K_c > 0$. We can now calculate
\begin{align}
\label{xmodel_hyp}
\partial_t d_2 &= \frac{1}{2} Q_f \partial_x \left(\frac{K_c(d_1,d_2)}{K_f(d_1,d_2)} \right)
\end{align}
The unknown $d_2$ gets transported with flux $Q_f K_c / K_f$ and can show hyperbolic behavior, such as the formation of discontinuities. 
\end{remark}

\subsection{Asymptotic Consistency}
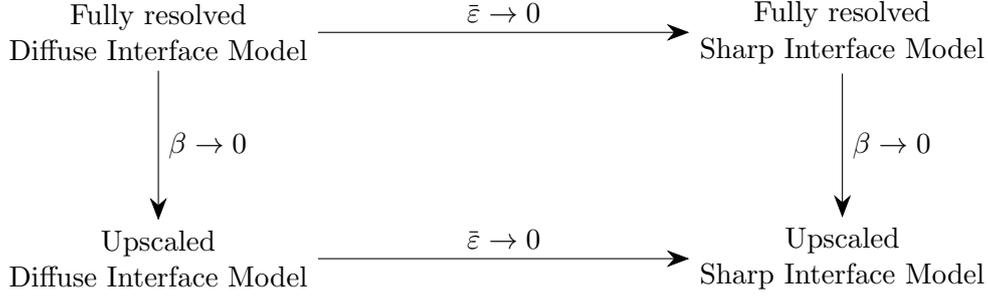
\begin{figure}
\centering
\begin{tikzpicture}
\node[align=center](n1){Fully resolved \\ Diffuse Interface Model};
\node[align=center, below of=n1, node distance = 3cm](n3){Upscaled \\ Diffuse Interface Model};
\draw[-{Stealth[length=3mm]}] (n1.south) --
node[right, align=left]{$\beta\to 0$} 
 (n3.north);
\node[align=center, right of=n1, node distance = 9cm](n2){Fully resolved \\ Sharp Interface Model};
\draw[-{Stealth[length=3mm]}] (n1.east) --
node[above, align=left]{$\overbar\eps\to 0$} 
 (n2.west);
\node[align=center, below of=n2, node distance = 3cm](n4){Upscaled \\ Sharp Interface Model};
\draw[-{Stealth[length=3mm]}] (n3.east) --
node[above, align=left]{$\overbar\eps\to 0$} 
 (n4.west);
\draw[-{Stealth[length=3mm]}] (n2.south) --
node[right, align=left]{$\beta\to 0$} 
 (n4.north);
\end{tikzpicture}
\caption{Models obtained by upscaling ($\beta \to 0$) and sharp interface limit ($\overbar \eps \to 0$).}
\label{fig:asy}
\end{figure}

In Section \ref{sec:up} we have investigated the limit process $\beta \to 0$, while in Section \ref{sec:sil} we examined $\overbar \eps \to 0$. An common question is under which circumstances there is asymptotic consistency, i.e. these two limit processes commute. In figure \ref{fig:asy} all limit processes are shown in a commutative diagram.

We investigate asymptotic consistency with non-dimensional numbers chosen as in \eqref{eq_nr_re}-\eqref{eq_nr_pec} with $\overbar{Re}$, $\overbar{Ca}$, $\overbar{Da}$, $\overbar M$, $\overbar{Pe_c}$ constant and independent of $\overbar \eps$ and $\beta$. The non-dimensional $\delta$ is chosen as $\delta = \overbar \eps$.

When starting with the fully-resolved diffuse-interface model \eqref{nmodel1}-\eqref{nmodel5} the limit $\overbar \eps \to 0$ results in a sharp-interface model as described in Section \ref{sec:si}. For details on this sharp-interface limit, see \cite{RvW}.

When we assume the geometry of Section \ref{sec:us_geo} we can proceed to upscale the fully-resolved sharp-interface model after introducing $d_1$ and $d_2$. While the process is tedious, the main ideas are analog to the calculations in \cite{Sohely}. In particular the asymptotic expansion of interface conditions, normal vectors and curvature has to be handled with care, as the coordinates $\xv = (x,\beta y)$ depend on $\beta$. For sake of brevity we skip this calculation here.

With the geometry of Section \ref{sec:us_geo} we find asymptotic consistency, that is the limit processes $\beta \to 0$ and $\overbar \eps \to 0$ commute. The result of the upscaling of the fully-resolved sharp-interface model is exactly given by \eqref{xmodel1}-\eqref{xmodel5}.

\begin{remark}
In more general geometries, asymptotic consistency does not necessary hold. This is due to the following observation. When upscaling a fully-resolved diffuse-interface model, the parameter $\delta$ is constant and independent of $\beta$. This leads to $\tilde \phi_f > 0$ and $\tilde \phi_c > 0$ everywhere. Because of this, we obtain upscaled equations for $p$ and $c$ without further assumptions on the geometry. The upscaled variables $p$ and $c$ do not depend on $y$, even if the geometry consists of two parallel channels separated by a solid region with $\Phiv \approx \ev_3$.
On the other hand, when upscaling the fully-resolved diffuse-interface model, the $\delta$-modifications have already vanished, as $\delta = \overbar \eps$. In this case, it is possible to have a different pressure $p$ in each channel, that is in each connected part of $\overline{\Omega_1|_x \cup \Omega_2|_x}$. Also it is possible to have a different ion concentration $c$ in each connected part of $\Omega_1|_x$.

We conclude that we have asymptotic consistency under the condition that there is only one flow channel, i.e. $\overline{\Omega_1|_x \cup \Omega_2|_x}$ is connected for every $x$, and that the first fluid phase is connected, i.e. $\Omega_1|_x$ is connected for every x. It is also possible to consider the symmetric case as in Section \ref{sec:us_geo} and have two symmetric connected parts of fluid one.
\end{remark}

\section{Numerical Investigation}
\label{sec:num}

We will now compare the upscaled $\delta$-$2f1s$-model \eqref{umodel1}-\eqref{umodel_v2} to the fully-resolved $\delta$-$2f1s$-model \eqref{nmodel1}-\eqref{nmodel5}. Remark \ref{rem:hyp} suggests that shock fronts can form in the upscaled model. Note that in this case the assumptions for the upscaling in Section \ref{sec:up} are no longer valid, and we expect different behaviours from the two models.

For the fully-resolved $\delta$-$2f1s$-model \eqref{nmodel1}-\eqref{nmodel5} we use a monolithic finite-element implementation provided by the DUNE-Phasefield module \cite{darusLars}. We employ Taylor--Hood elements for the flow variables velocity and pressure, and first-order Lagrange elements for the ion concentration and the phase-field parameters. The implementation is based on DUNE-PDELab \cite{pdelab} using ALU-Grid routines for adaptive grid generation \cite{alugrid}.

\subsection{Numerical Scheme for the Upscaled $\delta$-$2f1s$-Model}
\label{sec:num_up}
The upscaled $\delta$-$2f1s$-model consists of multiple coupled problems. The upscaled equations \eqref{umodel1}-\eqref{umodel3} for the unknowns $Q_f$, $p$ and $c$ have parameters \eqref{umodel_phictot}-\eqref{umodel_rtot} that depend on the distribution of phases in $y$-direction. This distribution is described by the fully coupled 2-d problem \eqref{umodel4}-\eqref{umodel5} for the Cahn--Hilliard variables $\phi_1$, $\phi_2$, $\phi_3$, $\mu_1$, $\mu_2$, $\mu_3$. Furthermore the flow profile has to be calculated by the cell problem \eqref{umodel2},\eqref{umodel2_bc}.

For simplicity we present the numerical scheme for equidistant time steps $t_n = n \Delta t$ and equidistant discretization in $x$ by $x_k = k \Delta x$. Let also $x_{k+1/2} = (x_k + x_{k+1})/2$. For each $t_n$, $x_k$ we discretize the one-dimensional unknown $\phi^{n}_{1,k}(y) = \phi_1(t_n,x_k,y)$ with linear Lagrange elements, and analogous for $\phi^{n}_{2,k}$, $\phi^{n}_{3,k}$, $\mu^{n}_{1,k}$, $\mu^{n}_{2,k}$, $\mu^{n}_{3,k}$, $\vv^{(1),n}_{0,k}$, $\vv^{(2),n}_{1,k}$, $w^n_k$. Again, we also use this notation for other variables such as $\tilde \phi^n_{f,k}$.

We discretize the macroscopic unknown $c(t,x)$ with a finite volume scheme, that is $c^n(x) = c(t_n,x)$ is piecewise constant with $c(t_n,x) = c^n_k$ for $x\in(x_{k-1/2},x_{k+1/2})$. The pressure $p^n(x) = p(t_n,x)$ is discretized using linear Lagrange elements with nodes $x_{k+1/2}$. Therefore $\partial_x p$ is constant on each finite volume cell $(x_{k-1/2},x_{k+1/2})$.

Given $\Phiv^n_k$, $c^n_k$ for all $x_k$ at time $t_n$, we now calculate the next time step using the following algorithm.
\begin{enumerate}
\item \label{step_w} For each $x_k$ use \eqref{umodel2},\eqref{umodel2_bc} to solve for $w^n_k(y)$. Here we use $\Phiv = \Phiv^n_k$ and the finite element method to discretize the equation. The equations for different $x_k$ are independent and can be solved in parallel.
\item For each $x_k$ calculate $K^n_{f,k}$ and $K^n_{c,k}$ by
\begin{align*}
K^n_{f,k} = \int_{-\ell_\Omega/2}^{\ell_\Omega/2} \tilde \phi^n_{f,k} w^n_k \diff y, \qquad K^n_{c,k} = \int_{-\ell_\Omega/2}^{\ell_\Omega/2} \tilde \phi^n_{c,k} w^n_k \diff y
\end{align*}
\item \label{step_K} Solve for $p^n(x)$ using the finite element method with
\begin{align*}
\partial_x( - K^n_f \partial_x p^n) = 0
\end{align*}
Here $K^n_f(x) = K^n_{f,k}$ for $x\in(x_{k-1/2},x_{k+1/2})$. As $K^n_f > 0$, the pressure $p$ is either a monotone increasing or monotone decreasing function, depending on the boundary conditions. We assume from here on $\partial_x p^n \leq 0$ and therefore fluid flow in positive $x$ direction. In case $\partial_x p^n \geq 0$ the upwind schemes in Steps \ref{step_pf} and \ref{step_c} have to be modified.
\item For each $x_k$ calculate $\vv^{(1),n}_{0,k}(y) = - w^n_k(y) \partial_x p^n(x_k)$.
\item \label{step_pf} Next, for each $x_k$ we solve for $\vv^{(2),n}_{1,k}$ and the Cahn-Hilliard variables $\phi^{n+1}_{2,k}$, $\phi^{n+1}_{3,k}$, $\mu^{n+1}_{1,k}$, $\mu^{n+1}_{2,k}$, $\mu^{n+1}_{3,k}$. For $\vv^{(2),n}_{1,k}$ we use \eqref{umodel_v2} the with an explicit upwind scheme for the $x$-derivative, i.e.,
\begin{align}
\label{nummodel_v2}
\partial_y (\tilde \phi^{n+1}_{f,k} \vv^{(2),n}_{1,k}) = - \frac{\tilde \phi^n_{f,k} \vv^{(1),n}_{0,k} - \tilde \phi^n_{f,k-1} \vv^{(1),n}_{0,k-1} }{\Delta x}.
\end{align}
This equation is coupled with the Cahn--Hilliard cell problems \eqref{umodel4}-\eqref{umodel5}. We again use an explicit upwinding scheme for the $x$-derivative, that is
\begin{align}
\begin{split}
&\frac{\phi^{n+1}_{1,k}-\phi^{n}_{1,k}}{\Delta t} + 
\frac{\phi^n_{1,k} \vv^{(1),n}_{0,k} - \phi^n_{1,k-1} \vv^{(1),n}_{0,k-1}}{\Delta x} +
\partial_y  (\phi^{n+1}_{1,k} \vv^{(2),n}_{1,k}) - \frac{\overbar \eps \overbar M}{\Sigma_1}  \partial_y^2 \mu^{n+1}_{1,k}  \\
&\qquad= - \frac{\overbar{Da}}{\overbar \eps } q(\Phiv^{n+1}_k)\left(r(c^n(x_k)) + \tilde \alpha \mu^{n+1}_{1,k} - \tilde \alpha \mu^{n+1}_{3,k}\right)
\end{split} \\
\begin{split}
&\frac{\phi^{n+1}_{2,k}-\phi^{n}_{2,k}}{\Delta t} + 
\frac{\phi^n_{2,k} \vv^{(1),n}_{0,k} - \phi^n_{2,k-1} \vv^{(1),n}_{0,k-1}}{\Delta x} +
\partial_y  (\phi^{n+1}_{2,k} \vv^{(2),n}_{1,k}) - \frac{\overbar \eps \overbar M}{\Sigma_1}  \partial_y^2 \mu^{n+1}_{2,k} = 0
\end{split} \\
\label{nummodel_phi3}
&\phi^{n+1}_{3,k} = 1- \phi^{n+1}_{1,k} -\phi^{n+1}_{2,k} \\
&\mu^{n+1}_{1,k} = \frac{\partial_{\phi_1} W(\Phiv^{n+1}_k)}{\overbar \eps} - \overbar \eps \Sigma_i \partial_y^2 \phi^{n+1}_{1,k} \\
&\mu^{n+1}_{2,k} = \frac{\partial_{\phi_2} W(\Phiv^{n+1}_k)}{\overbar \eps} - \overbar \eps \Sigma_i \partial_y^2 \phi^{n+1}_{2,k} \\
\label{nummodel_mu3}
&\mu^{n+1}_{3,k} = - \mu^{n+1}_{1,k} - \mu^{n+1}_{2,k}
\end{align}
Note that we do not use \eqref{umodel4b} and \eqref{umodel5} for $\phi^{n+1}_{3,k}$ and $\mu^{n+1}_{3,k}$. Instead we use that by construction $\phi_1 + \phi_2 + \phi_3 = 1$ and $\mu_1 + \mu_2 + \mu_3 = 0$, see \cite{RvW} for details.

We use the finite element method to discretize \eqref{nummodel_v2}-\eqref{nummodel_mu3} and Newtons method to solve the resulting nonlinear system. This step has by far the highest computational cost. With the explicit upwinding scheme for the $x$ derivatives, the cell problems for each $k$ fully decouple and can be solved in parallel. This leads to a significant speedup in comparison to discretizing the Cahn--Hilliard evolution \eqref{umodel4}-\eqref{umodel5} naively as a 2-d problem. 
\item \label{step_R} Calculate $\tilde \phi^{n+1}_{c,\text{total},k}$ and $R^{n+1}_{\text{total},k}$ as
\begin{align}
\tilde \phi^{n+1}_{c,\text{total},k} &= \int_{-\ell_\Omega/2}^{\ell_\Omega/2} \phi^{n+1}_{c,k} \diff y \\
R^{n+1}_{\text{total},k} &= - \int_{-\ell_\Omega/2}^{\ell_\Omega/2} q(\Phiv^{n+1}_k)\left(r(c^n(x_k)) + \tilde \alpha \mu^{n+1}_{1,k} - \tilde \alpha \mu^{n+1}_{3,k}\right)\diff y
\end{align}
We also set $\tilde \phi^{n+1}_{c,\text{total},k+1/2} = (\tilde \phi^{n+1}_{c,\text{total},k} +\tilde \phi^{n+1}_{c,\text{total},k+1})/2$.
\item \label{step_c} Finally we solve for $c$ using \eqref{umodel3} discretized by the finite volume method. We use an implicit upwinding scheme for the transport in $x$-direction
\begin{align}
\begin{split}
&\frac{\tilde \phi^{n+1}_{c,\text{total},k} c^{n+1}_k - \tilde \phi^{n}_{c,\text{total},k} c^{n}_k}{\Delta t} 
- \frac{ K^{n}_{c,k} \partial_x p^n(x_k) c^{n+1}_k -  K^{n}_{c,k-1} \partial_x p^n(x_{k-1}) c^{n+1}_{k-1}}{\Delta x} \\
&\qquad = \frac{1}{\overbar{Pe_c}} \frac{1}{\Delta x} \left( \tilde \phi^{n+1}_{c,\text{total},k+1/2} \frac{c^{n+1}_{k+1} - c^{n+1}_{k}}{\Delta x} - \phi^{n+1}_{c,\text{total},k-1/2} \frac{c^{n+1}_k - c^{n+1}_{k-1}}{\Delta x}\right)
+ \frac{\overbar{Da}}{\overbar \eps } R^{n+1}_{\text{total},k}
\end{split}
\end{align}
\end{enumerate}

\subsection{Comparison: Formation of an $N$-Wave}
\label{sec:num_shock}
As our first numerical example we choose a geometry as described in Section \ref{sec:us_geo}, with the computational domain $(x,y) \in [0,1]\times[-1,0]$.  For $x = 0$ and $x = 1$ we choose periodic boundary conditions for all variables except the pressure $p$. For $y = -1$ we use the trivially upscaled versions of the boundary conditions \ref{bc_phi}-\ref{bc_v} and for $y=0$ we choose boundary conditions according to the symmetry assumption.

We will compare the non-dimensional $\delta$-$2f1s$-model with the upscaled $\delta$-$2f1s$-model \eqref{umodel1}-\eqref{umodel_v2}. For simplicity we choose $\gamma_1 = \gamma_2$ and $d_0$ sufficiently big such that $L_\text{slip} \approx 0$. We choose the phase-field parameter $\overbar \eps = 0.03$ and $\delta = \overbar \eps$ as in Section \ref{sec:sil}.

\begin{figure}
\includegraphics[height = 3.8cm, trim=50 0 600 0, clip]{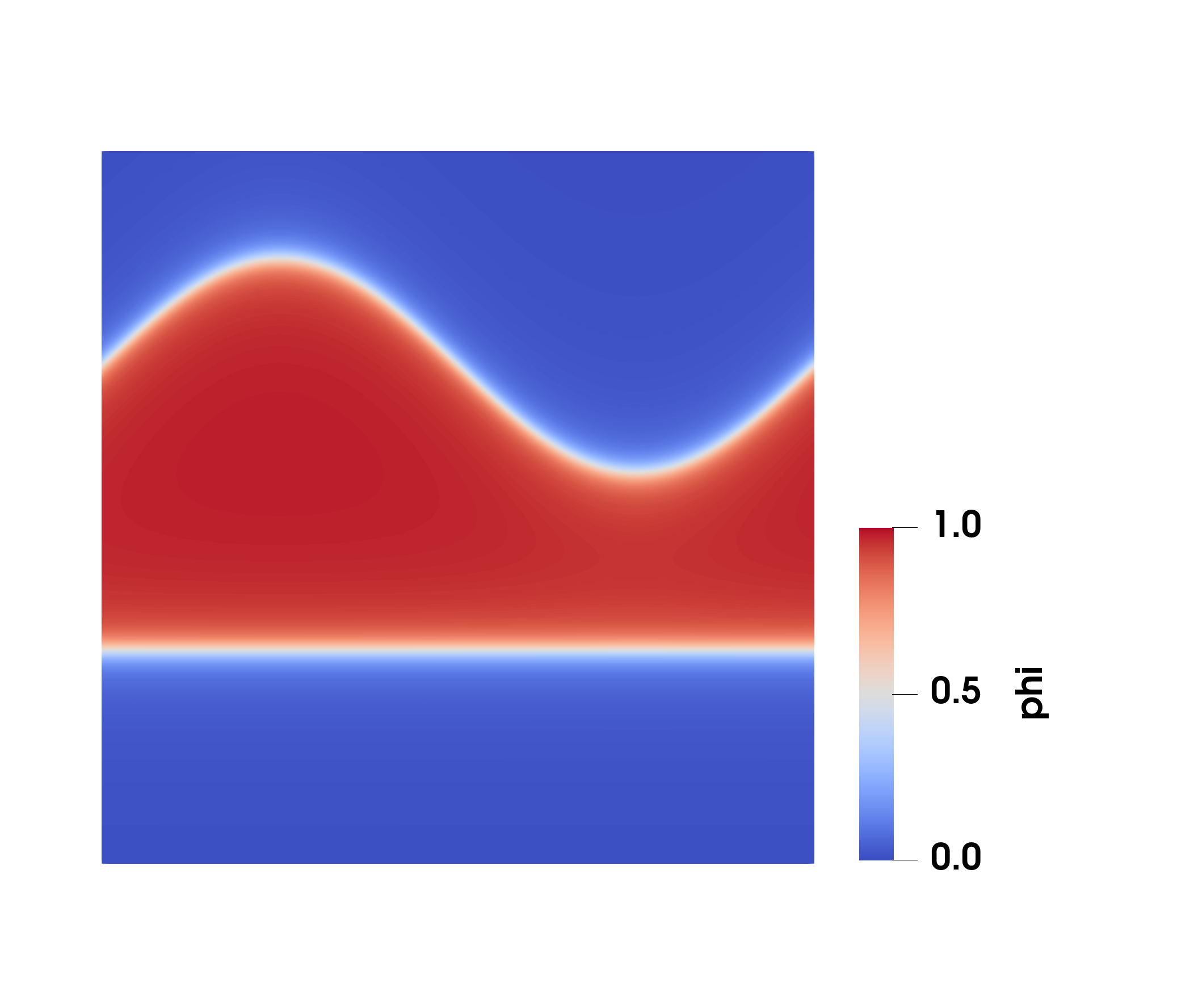}
\includegraphics[height = 3.8cm, trim=50 0 600 0, clip]{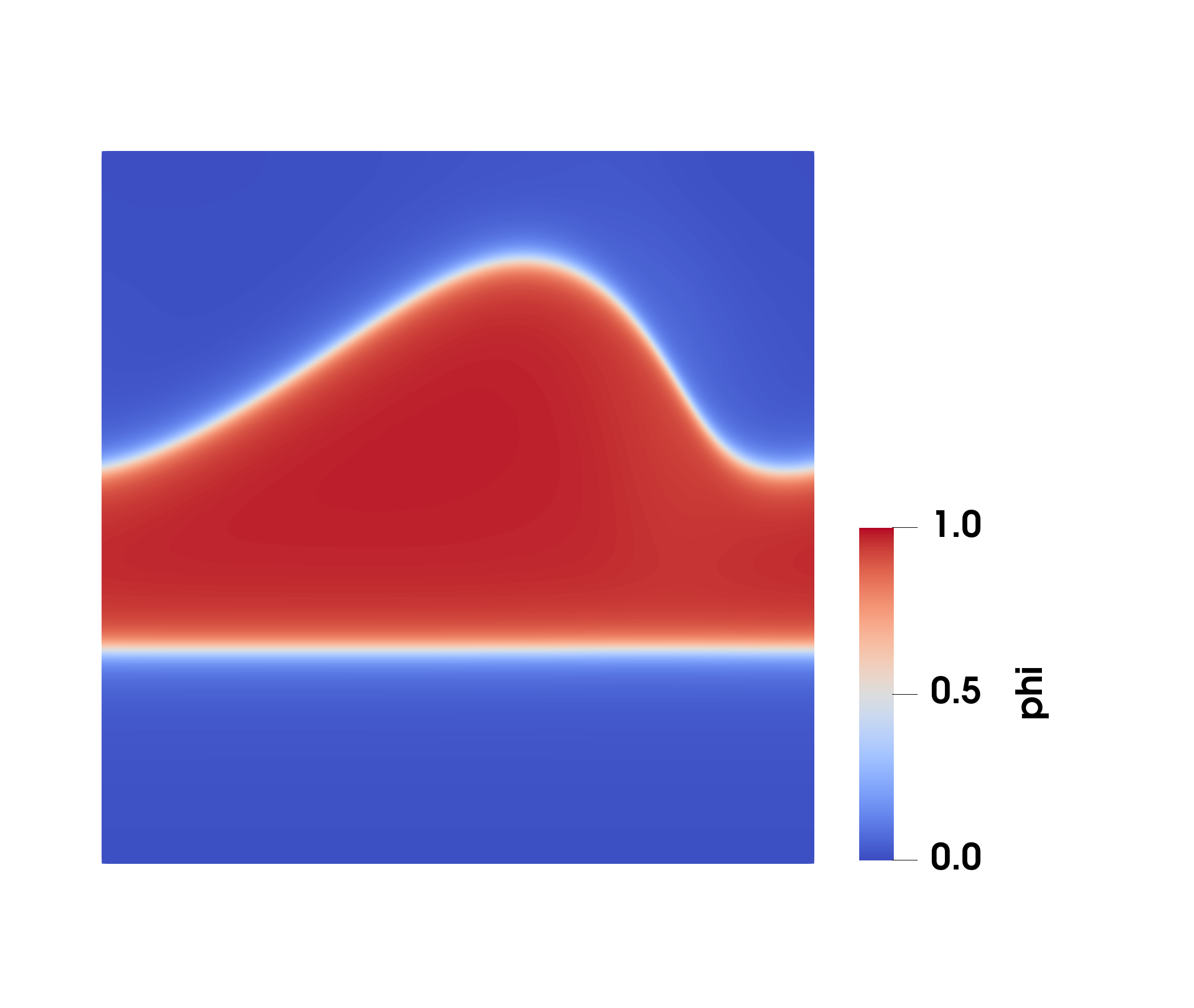}
\includegraphics[height = 3.8cm, trim=50 0 600 0, clip]{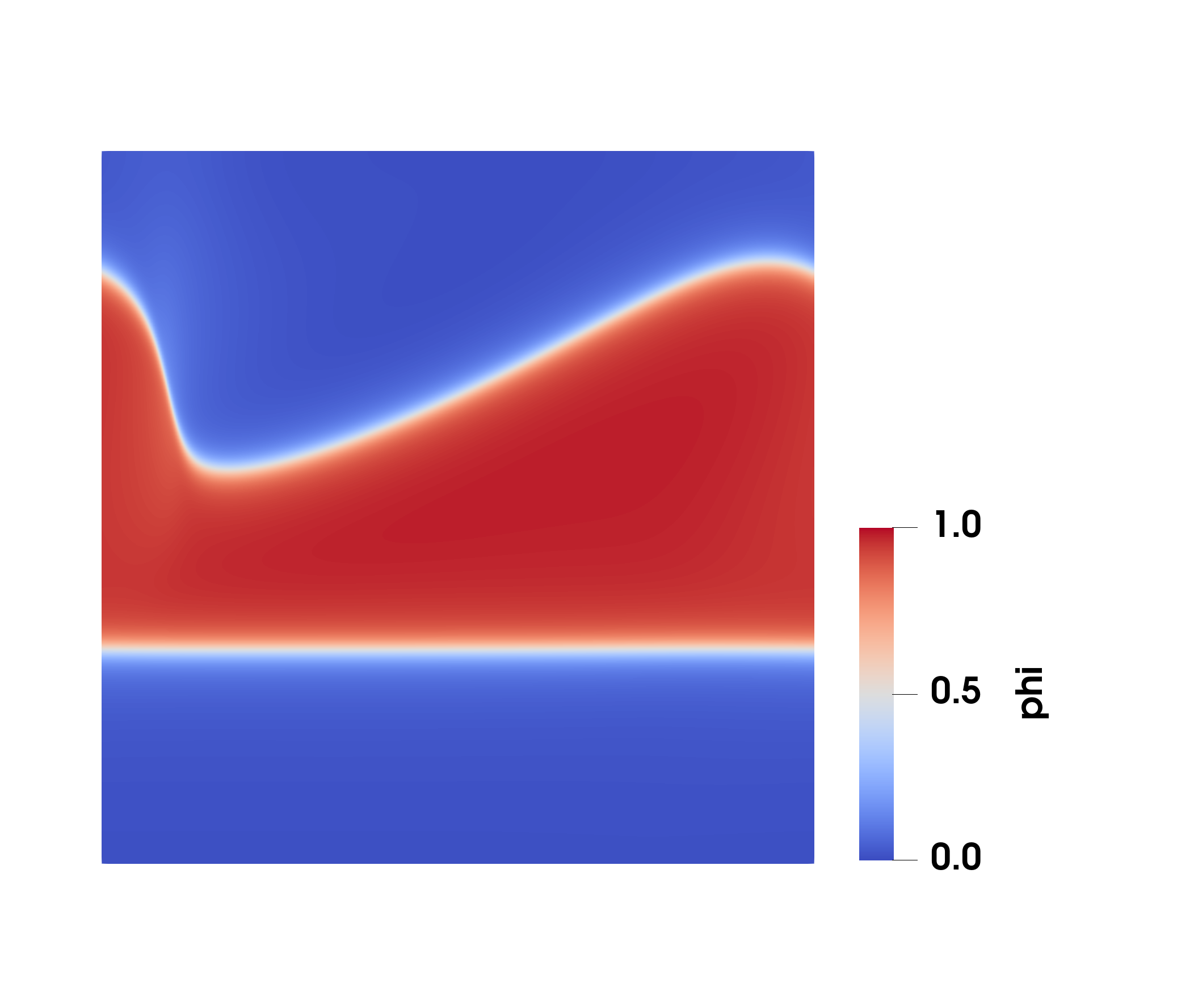}
\includegraphics[height = 3.8cm, trim=50 0 200 0, clip]{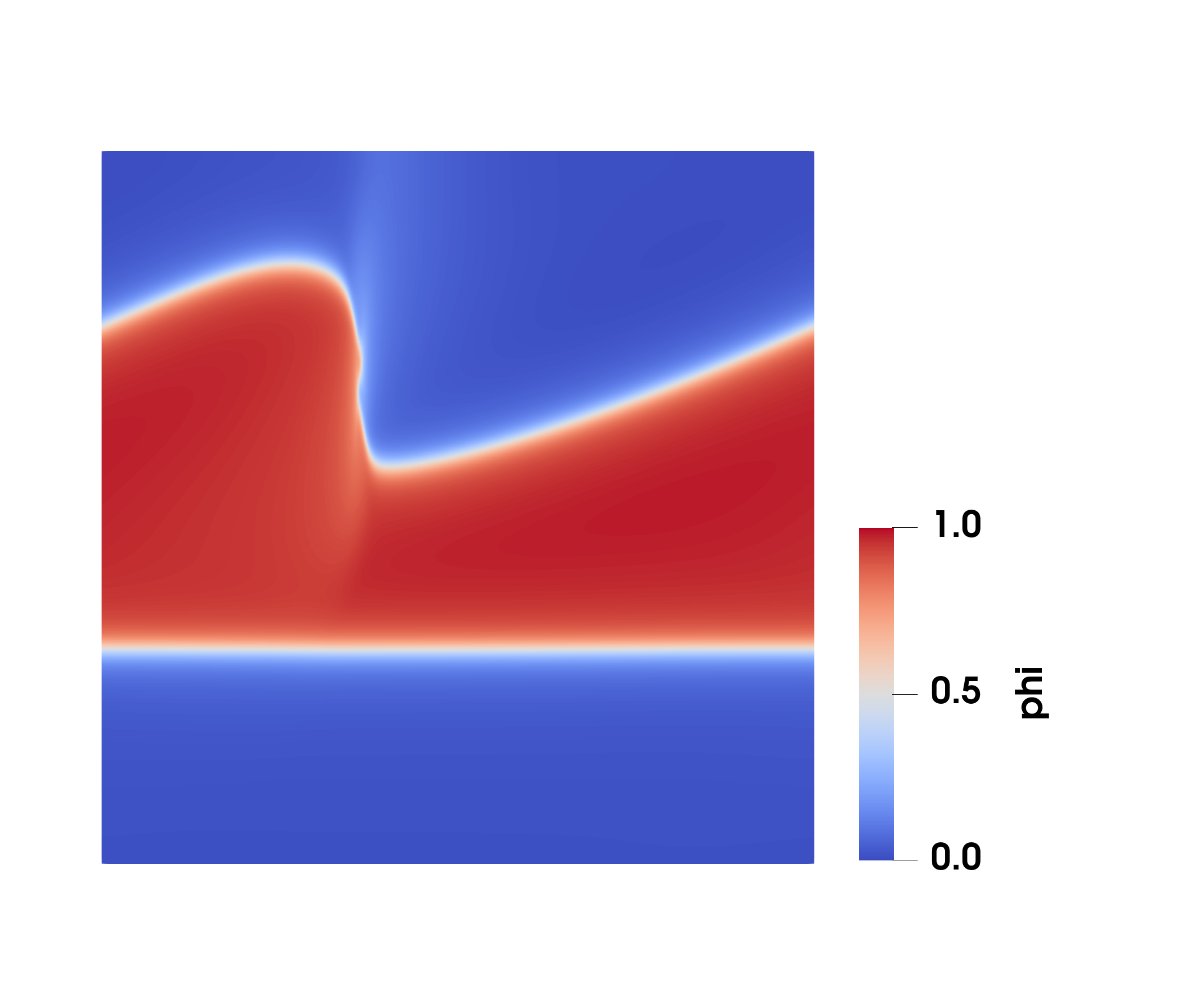}
\caption{Evolution of the upscaled $\delta$-$2f1s$-model on the domain $[0,1]\times[-1,0]$. Shown in red is fluid phase one, with fluid phase two above and solid phase below. From left to right: Initial data, $t=0.15$, $t=0.3$ and $t=0.45$.}
\label{fig:r1pic}
\end{figure}

We want to focus on the hyperbolic behavior of $d_2$ as described in Remark \ref{rem:hyp}. Therefore we choose $c$ in the initial conditions such that $r(c) = 0$. This leads to no precipitation or dissolution in the model, and the fluid-solid interface does not change over time. We choose 
\begin{align*}
d_1 + d_2 \equiv 0.7\qquad \text{and}\qquad d_1(x) = 0.4 + 0.15\sin(2\pi x).
\end{align*}
This corresponds to a plane fluid-solid interface and a sine-shaped fluid-fluid interface. An image of these initial conditions is given in figure \ref{fig:r1pic}. 

By applying a pressure difference as Dirichlet boundary condition at $x=0$ and $x=1$, the two fluid phases will move in positive $x$-direction. The fluid velocity $\vv^{(1)}_0$ is higher in the center of the channel. As shown in figure \ref{fig:r1pic} this will lead to a steeper fluid-fluid interface over time. At a time $t^\ast>0$ the upscaled $\delta$-$2fs$ model has a fluid-fluid interface that is perpendicular to the thin strip. As discussed in Remark \ref{rem:scaling}, the assumptions for the upscaling in Section \ref{sec:up} are no longer valid. For times $t>t^\ast$ the fluid-fluid interface will roll over, leading to multiple layers of fluid phase 1 at the same $x$ value.

\begin{figure}
\centering
\begin{tikzpicture}
\begin{axis}[width=6.2 cm, height = 6.2 cm, xmin=0, xmax=1, no markers,
xlabel={$x$},ylabel={$y$},
every axis y label/.style={
    at={(ticklabel* cs:1.0)},
    anchor=east,
}]
\addplot[solid, 			color=red!100!black] table 
[col sep=comma, x index= 0, y expr={-1+\thisrowno{1}}] {data1.csv};
\addplot[densely dashed, color=red!80!black] table 
[col sep=comma, x index= 0, y expr={-1+\thisrowno{2}}] {data1.csv};
\addplot[dashed,			color=red!60!black] table 
[col sep=comma, x index= 0, y expr={-1+\thisrowno{3}}] {data1.csv};
\addplot[dashdotted,		color=red!40!black] table 
[col sep=comma, x index= 0, y expr={-1+\thisrowno{4}}] {data1.csv};
\addplot[densely dotted,	color=red!20!black] table 
[col sep=comma, x index= 0, y expr={-1+\thisrowno{5}}] {data1.csv};
\addplot[solid] table 
[col sep=comma, x index= 0, y expr={-1+\thisrowno{6}}] {data1.csv};
\end{axis}
\end{tikzpicture} \hspace{0.0cm}
\begin{tikzpicture}
\begin{axis}[width=6.2 cm, height = 6.2 cm, xmin=0, xmax=1, no markers,
legend style={cells={anchor=east},legend pos=outer north east},
xlabel={$x$},ylabel={$y$},
every axis y label/.style={
    at={(ticklabel* cs:1.0)},
    anchor=east,
}]
\addplot[solid, 			color=red!100!black] table 
[col sep=comma, x index= 0, y expr={-1+\thisrowno{1}}] {data2.csv};
\addplot[densely dashed, color=red!80!black] table 
[col sep=comma, x index= 0, y expr={-1+\thisrowno{2}}] {data2.csv};
\addplot[dashed,			color=red!60!black] table 
[col sep=comma, x index= 0, y expr={-1+\thisrowno{3}}] {data2.csv};
\addplot[dashdotted,		color=red!40!black] table 
[col sep=comma, x index= 0, y expr={-1+\thisrowno{4}}] {data2.csv};
\addplot[densely dotted,	color=red!20!black] table 
[col sep=comma, x index= 0, y expr={-1+\thisrowno{5}}] {data2.csv};
\addplot[solid] table 
[col sep=comma, x index= 0, y expr={-1+\thisrowno{6}}] {data2.csv};
\legend{$\beta=1$,$\beta=\frac{1}{2}$,$\beta=\frac{1}{4}$,$\beta=\frac{1}{8}$,$\beta=\frac{1}{16}$,up.}
\end{axis}
\end{tikzpicture}
\caption{Fluid-fluid interface locations for the non-dimensional $\delta$-$2f1s$-model with varying $\beta$, and for the upscaled $\delta$-$2f1s$-model. The interface is located through the condition $\phi_1 = \phi_2$. Left: $t=0.3$, Right: $t=0.44$.}
\label{fig:r1comparison}
\end{figure}
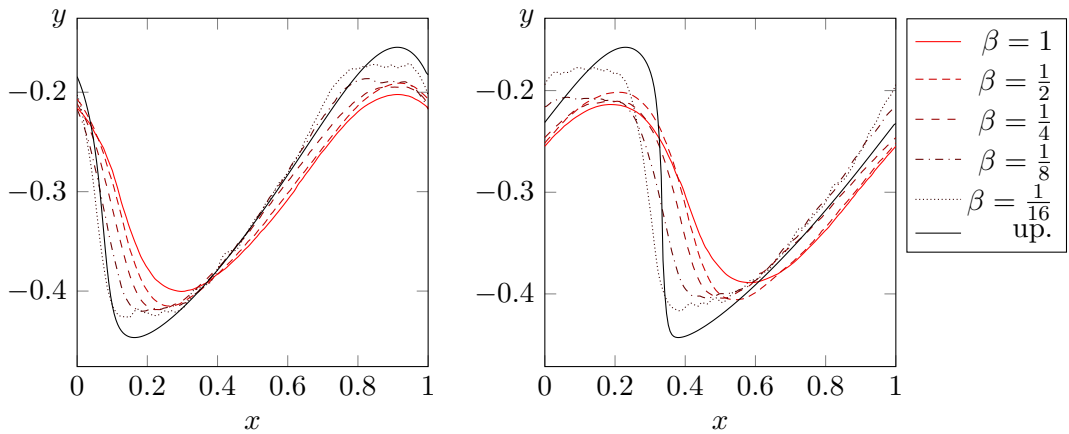

We can compare this behavior with the non-dimensional $\delta$-$2f1s$-model in a thin strip for different values of $\beta$. As shown in figure \ref{fig:r1comparison}, for times $t < t^\ast$ there is a good agreement between the non-dimensional $\delta$-$2f1s$-model with small values of $\beta$ and the upscaled $\delta$-$2f1s$-model. 

In contrast to the upscaled $\delta$-$2f1s$-model, the non-dimensional $\delta$-$2f1s$-model does not evolve to a fluid-fluid interface perpendicular to the thin strip, as shown in figure \ref{fig:r1comparison}. Instead, 
when reaching a steep fluid-fluid interface there are regions of high curvature at the beginning and end of the steep passage. In these regions of high curvature the surface tension leads to a pressure difference between the fluid phases, counteracting the interface getting steeper. For smaller $\beta$ the fluid-fluid interface allows for a steeper passage in $(x,y)$ coordinates, as this effect depends on the curvature in the $\xv$ coordinates, which are not scaled with $\beta$.

\subsection{Comparison: Precipitation}
In the second numerical example we study precipitation in the thin strip. We use the same domain and boundary conditions as in the previous example. Again, we choose $\gamma_1 = \gamma_2$, and a $d_0$ large enough so that $L_\text{slip}\approx 0$. We further choose $\overbar \eps = 0.03$ and $\delta = \overbar \eps$.
We use a simple, linear reaction rate $r(c) = c - 0.5$ and choose the ion concentration to be in equilibrium initially, that is $c = 0.5$ everywhere. With $d_1(x) = 0.4$ and $d_2(x) = 0.3$ in the initial conditions correspond to the phases  being layered in the thin strip, without depending on $x$.
To induce precipitation we add a source term $s(x)$ to the ion conservation equation \eqref{nmodel3}, it now reads
\begin{align*}
 \partial_{t} (\tilde\phi_c c) + \nabla \cdot (\phi_c \vv c) + \frac{Cn}{\beta Pe_{CH}} \nabla \cdot(\Jv_1 c) &= \frac{1}{Pe_c} \nabla \cdot ( \tilde\phi_c  \nabla c) + Da R + \tilde \phi_c s(x).
\end{align*}
The source terms upscales trivially at order $O(\beta^0)$, and the upscaled ion conservation equation \eqref{umodel3} is now given by
\begin{align*}
\frac{d}{dt} \left( \tilde \phi_{c,\text{total}} c \right) + \partial_x \left( (-K_c \partial_x p) c \right) &= \frac{1}{\overbar{Pe_c}} \partial_x \left( \tilde \phi_{c,\text{total}} \partial_x c\right) + \frac{\overbar{Da}}{\overbar \eps } R_\text{total} 
+ \tilde \phi_{c,\text{total}} s(x).
\end{align*}
We choose the ion source to be located between $x = 0.1$ and $x=0.3$, in detail
\begin{align*}
s(x) = \max(0, 62.5 (x-0.1)(0.3-x))
\end{align*}

\begin{figure}
\centering
\begin{tikzpicture}
\begin{axis}[width=7.0 cm, height = 7.0 cm, xmin=0, xmax=1, no markers,
legend style={cells={anchor=east},legend pos=outer north east},
xlabel={$x$},ylabel={$y$},
every axis y label/.style={
    at={(ticklabel* cs:1.0)},
    anchor=east,
}]
\addplot[solid, 			color=red!100!black] table 
[col sep=comma, x index= 0, y expr={-1+\thisrowno{1}}] {data3_1.csv};
\addplot[densely dashed, color=red!80!black] table 
[col sep=comma, x index= 0, y expr={-1+\thisrowno{2}}] {data3_1.csv};
\addplot[dashed,			color=red!60!black] table 
[col sep=comma, x index= 0, y expr={-1+\thisrowno{3}}] {data3_1.csv};
\addplot[dashdotted,		color=red!40!black] table 
[col sep=comma, x index= 0, y expr={-1+\thisrowno{4}}] {data3_1.csv};
\addplot[densely dotted,	color=red!20!black] table 
[col sep=comma, x index= 0, y expr={-1+\thisrowno{5}}] {data3_1.csv};
\addplot[solid, 			color=red!100!black] table 
[col sep=comma, x index= 0, y expr={-1+\thisrowno{1}}] {data3_2.csv};
\addplot[densely dashed, color=red!80!black] table 
[col sep=comma, x index= 0, y expr={-1+\thisrowno{2}}] {data3_2.csv};
\addplot[dashed,			color=red!60!black] table 
[col sep=comma, x index= 0, y expr={-1+\thisrowno{3}}] {data3_2.csv};
\addplot[dashdotted,		color=red!40!black] table 
[col sep=comma, x index= 0, y expr={-1+\thisrowno{4}}] {data3_2.csv};
\addplot[densely dotted,	color=red!20!black] table 
[col sep=comma, x index= 0, y expr={-1+\thisrowno{5}}] {data3_2.csv};
\legend{$\beta=1$,$\beta=\frac{1}{2}$,$\beta=\frac{1}{4}$,$\beta=\frac{1}{8}$,up.}
\end{axis}
\end{tikzpicture}
\caption{Interface locations at time $t=2.4$ for the non-dimensional $\delta$-$2f1s$-model with varying $\beta$, and for the upscaled $\delta$-$2f1s$-model. The fluid-fluid interface can be seen in the upper half and is located by the condition $\phi_1 = \phi_2$. The fluid-solid interface in the lower half is located by $\phi_1 = \phi_3$.}
\label{fig:r2comparison}
\end{figure}
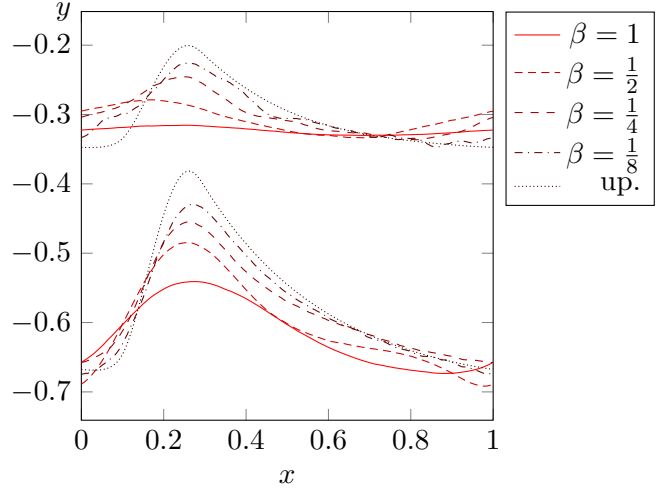

Figure \ref{fig:r2comparison} shows a comparison between the non-dimensional $\delta$-$2f1s$-model with different values of $\beta$, and the upscaled $\delta$-$2f1s$-model. There is a good agreement between the full model with small values of $\beta$ and the upscaled model. For large values of $\beta$ there is less precipitation in the thin strip. This is due to the ion concentration $c$ not being constant in $y$-direction. The source term $\tilde \phi_c s(x)$ generates ions everywhere in the first fluid phase, but precipitation removes ions from the first fluid phase only at the fluid-solid interface. This leads to an oversaturation $c>0.5$ further away from the fluid-solid interface. For smaller values of $\beta$ the diffusion in $y$-direction results in more ions precipitating and therefore a smaller oversaturation of ions in the fluid phase.

Figure \ref{fig:r2comparison} also shows the influence of a non-constant width of the thin strip on the flow inside the thin strip. The fluid-fluid interfaces are pushed towards the center of the thin strip, where flow velocities are higher.

\bibliographystyle{siamplain}
\bibliography{references}
\end{document}